\documentclass[11pt,a4paper]{amsart}

\usepackage{amsmath,amssymb,amsthm}
\usepackage{geometry}
\usepackage{enumitem}
\usepackage{mathtools}
\usepackage{float}

\usepackage[no-math]{fontspec}

\usepackage{tikz}
\usetikzlibrary{calc, arrows.meta, positioning, decorations.pathreplacing}

\usepackage{xurl}
\usepackage[colorlinks=true, linkcolor=blue, citecolor=red, urlcolor=blue]{hyperref}
\tikzset{
    grid line/.style={step=1cm, gray!30, very thin},
    nucleus/.style={fill=gray!10, draw=gray!50, ultra thick, dashed, rounded corners=2pt},
    toast piece/.style={fill=blue!5, draw=blue!60, ultra thick, rounded corners=4pt},
    sub piece/.style={fill=gray!20, draw=gray!60, ultra thick, rounded corners=4pt},
    node std/.style={circle, fill=black!70, inner sep=0pt, minimum size=6pt},
    node root/.style={rectangle, fill=purple!70, inner sep=0pt, minimum size=7pt, rounded corners=1pt},
    node inner/.style={circle, fill=red!70, inner sep=0pt, minimum size=6pt},
    node outer/.style={circle, fill=blue!70, inner sep=0pt, minimum size=6pt},
    node interior/.style={circle, fill=gray!50, inner sep=0pt, minimum size=6pt},
    path edge/.style={thick, black!30},
    path tree/.style={->, >=stealth, ultra thick, blue!80!black},
    path detour/.style={ultra thick, green!60!black, rounded corners=2pt},
    dist arrow/.style={<->, >=stealth, thick, red!80!black, dashed}
}

\newtheorem{theorem}{Theorem}[section]
\newtheorem{lemma}[theorem]{Lemma}
\newtheorem{claim}[theorem]{Claim}
\newtheorem{proposition}[theorem]{Proposition}
\newtheorem{corollary}[theorem]{Corollary}

\newtheorem{remark}[theorem]{Remark}
\newtheorem{question}[theorem]{Question}
\newtheorem{example}[theorem]{Example}

\title{The 1--2--3 Conjecture for Borel Combinatorics}

\author{Junhao Chen}
\address{School of Mathematical Sciences and LPMC, Nankai University, Tianjin 300071, P.R. China}
\email{jhchen@mail.nankai.edu.cn}

\author{Wei Dai}
\address{School of Mathematical Sciences and LPMC, Nankai University, Tianjin 300071, P.R. China}
\email{weidai@nankai.edu.cn}
\thanks{The second author acknowledges the partial support of their research by the Fundamental Research Funds for the Central Universities and by the National Natural Science Foundation of China (NSFC) grant 12271263.}

\date{}

\begin{document}

\begin{abstract}
The 1--2--3 Conjecture asserts that every finite graph without isolated edges admits edge weights from $\{1,2,3\}$ such that the sum function forms a vertex coloring. We study its Borel analogue for locally finite Borel graphs without isolated edges, proving that five weights suffice in general and four suffice for regular graphs. We obtain three weights for acyclic graphs admitting a Borel selection of a nonempty finite set of ends from each component. We also obtain $\mu$-measurable and Baire measurable vertex-coloring $3$-edge weightings from $3$-colorings of the same definability, under probability measure-preserving and Baire quasi-invariance conditions, respectively.
\end{abstract}

\maketitle

\section{Introduction}

Let $G=(X,E)$ be a simple undirected finite graph. A $k$\textit{-edge-weighting} is a symmetric function $w:E\longrightarrow\{1,\dots,k\}$. Given a $k$-edge-weighting $w$, we say that $w$ is a \textit{vertex-coloring} $k$-edge weighting if, for each vertex $x\in X$, the \textit{sum function} $s_w(x):=\sum_{y\in N_G(x)}w(xy)$ is a proper coloring of $G$. In 2004, Karo\'{n}ski, \L{}uczak, and Thomason conjectured that for every finite graph without isolated edges, there exists a vertex-coloring $3$-edge weighting~\cite{KLT04}. This is known as the 1--2--3 Conjecture. The exclusion of isolated edges is necessary, and three weights are needed in general: a triangle admits no vertex-coloring $2$-edge weighting.

This conjecture has received a lot of attention in the past two decades. Karo\'{n}ski, \L{}uczak, and Thomason verified the conjecture for $3$-colorable graphs in 2004~\cite{KLT04}. Three years later, Addario-Berry, Dalal, McDiarmid, Reed, and Thomason provided the first finite, general upper bound of $k=30$~\cite{ABDMRT07}. In 2010, Kalkowski, Karo\'{n}ski, and Pfender decreased the upper bound to $k=5$, using a technical algorithmic argument~\cite{KKP10}.  For $d$-regular graphs, a bound of $k=4$ was proved in 2021~\cite{Prz21} and later for general graphs in 2023~\cite{Keu23}. Finally, the conjecture was settled by Keusch in 2024.

\begin{theorem}[Keusch~\cite{Keu24}]
\label{keusch}
For every finite simple graph $G=(X,E)$ with no isolated edges, there is a weighting $w:E\longrightarrow\{1,2,3\}$ for which
\[
  \sum_{z\in N_G(x)} w(xz)
  \neq
  \sum_{z\in N_G(y)} w(yz)
  \qquad\text{whenever }xy\in E.
\]
\end{theorem}

Following Keusch's proof, Stawiski~\cite{Sta23} extended the result to locally finite infinite graphs without isolated edges, using a compactness argument. More recently, Charlesworth, Ramsey, and Strungaru~\cite{CRS26} gave explicit vertex-coloring $3$-edge weightings for several periodic and non-periodic polygonal tilings.  For a general overview of the progress on the 1--2--3 Conjecture and on related problems, we refer to the survey of Seamone~\cite{Sea12}.

This paper contributes to Borel combinatorics, which is a new branch of descriptive set theory. Here we give a brief introduction. Borel combinatorics mainly focuses on how graph theory concepts, like colorings and matchings, behave under definability constraints. It also has deep connections with other areas such as ergodic theory and probability theory. In 1999, Kechris, Solecki, and Todor\v{c}evi\'{c} published a celebrated paper on Borel chromatic numbers~\cite{KST99}, which is known as the birth of Borel combinatorics.  This area has developed over the last two decades; we refer readers to two surveys by Pikhurko~\cite{Pik21} and Kechris and Marks~\cite{KM20} for an overview of this new area.

In this paper, we introduce an analogue of the 1--2--3 Conjecture for Borel combinatorics. To ensure that the sum function is well-defined, we will assume that all graphs are locally finite. Suppose $G=(X,E)$ is a locally finite Borel graph. A \textit{Borel vertex-coloring $k$-edge weighting} is a vertex-coloring $k$-edge weighting which is additionally Borel. Since isolated vertices do not affect the edge-weighting problem, we may delete them and assume that the underlying graph has no isolated vertices. The following is the general question:

\begin{question}
For any locally finite Borel graph $G=(X,E)$ without isolated edges, does it admit a Borel vertex-coloring $3$-edge weighting?
\end{question}

Our first two results focus on finding a general upper bound for Borel graphs. Theorem~\ref{keusch} can be extended to the class of infinite graphs in the classical setting, but the proof is non-constructive and hence cannot be used in a Borel way. By a simple power graph argument, we show that there exists a general upper bound for all locally finite Borel graphs with no isolated edges.
\begin{theorem}
Every locally finite Borel graph with no isolated edges admits a Borel vertex-coloring $5$-edge weighting.
\end{theorem}

The proof adapts the algorithm of Kalkowski, Karo\'{n}ski, and Pfender~\cite{KKP10} to the Borel setting. The underlying method originates in Kalkowski's work on total weightings, presented in his doctoral thesis~\cite{Kal09}. There he constructed vertex-coloring total weightings with vertex weights in $\{1,2\}$ and edge weights in $\{1,2,3\}$. The central idea is to process vertices successively and make local changes so that previously established coloring constraints remain valid. Kalkowski, Karo\'{n}ski, and Pfender developed this method to obtain the bound $5$ for finite graphs~\cite{KKP10}. Another modification of Kalkowski's algorithm was used by Przyby\l{}o~\cite{Prz21} to obtain the bound $4$ for finite regular graphs. By adapting Przyby\l{}o's construction to the Borel setting, we also obtain the following result.

\begin{theorem}
Every regular locally finite Borel graph with no isolated edges admits a Borel vertex-coloring $4$-edge weighting.
\end{theorem}

Early approaches to the classical 1--2--3 Conjecture focused on the relationship between the number of edge weights and the chromatic number $\chi(G)$. A natural idea is to use a proper vertex coloring to guide the construction of an edge weighting, aiming to distinguish the color classes through the sums of incident edge weights. One way to achieve this is to arrange that vertices in different color classes have sums in different residue classes modulo a suitable integer. In the original paper of Karo\'{n}ski, \L{}uczak, and Thomason~\cite{KLT04}, they obtain the following result.

\begin{theorem}[Karo\'{n}ski, \L{}uczak, and Thomason~{\cite[Theorem~1]{KLT04}}]\label{kcolor}
Let $k$ be a positive odd integer. Every finite $k$-colorable graph with no isolated edges admits a vertex-coloring $k$-edge weighting.
\end{theorem}

It is natural to ask whether there is a similar relationship between the Borel chromatic number and Borel vertex-coloring edge weightings. This question motivates our subsequent results and their proofs. Under additional structural assumptions, we obtain Borel analogues of the classical result. As an application, we show that grid graphs arising from free Borel actions of $\mathbb{Z}^d$, for $d\geq2$, admit Borel vertex-coloring $3$-edge weightings.

\begin{theorem}
Suppose $G=(X, E)$ is a locally finite Borel graph with no isolated edges. If $G$ is Borel $k$-colorable, then $G$ has a Borel vertex-coloring $(k+1)$-edge weighting. Moreover, if $G$ admits a Borel one-ended spanning forest, then $G$ has a Borel vertex-coloring $k$-edge weighting.
\end{theorem}

The moreover part of the preceding theorem relies on one-ended spanning forests in constructing edge weightings. We next show that, for acyclic graphs, a Borel selection of a nonempty finite set of ends from each component suffices to obtain the bound $3$, without assuming Borel $3$-colorability.

\begin{theorem}
Suppose $G=(X, E)$ is a locally finite acyclic Borel graph with no isolated edges. If there is a Borel way to select a nonempty finite set of ends from each component of $G$, then there is a Borel vertex-coloring $3$-edge weighting.
\end{theorem}
We also consider measurable and Baire measurable edge weightings. Under the respective measure and category hypotheses, our final two results establish analogues of Theorem~\ref{kcolor} for $3$-colorable locally finite Borel graphs. First consider the probability measure-preserving (p.m.p.) setting.
\begin{theorem}
Suppose $G=(X,E)$ is a locally finite Borel graph that is probability measure-preserving with respect to a Borel probability measure $\mu$, and $G$ has no isolated edges. If $G$ has a $\mu$-measurable $3$-coloring, then $G$ admits a $\mu$-measurable vertex-coloring $3$-edge weighting.
\end{theorem}
A parallel result holds in the Baire category setting when the compatible Polish topology is $G$-quasi-invariant.

\begin{theorem}
Suppose $G=(X,E)$ is a locally finite Borel graph on a Polish space $(X,\tau)$ and $G$ has no isolated edges. Suppose that $\tau$ is $G$-quasi-invariant. If $G$ has a Baire measurable $3$-coloring, then $G$ admits a Baire measurable vertex-coloring $3$-edge weighting.
\end{theorem}

This paper is organized as follows. Section~\ref{sec:preliminaries} introduces the background in graph theory, descriptive set theory, and Borel combinatorics. In Section~\ref{sec:upper-bounds}, we establish the general bound $5$ for locally finite Borel graphs and improve it to $4$ for regular graphs. Section~\ref{sec:weightings} studies the role of Borel colorings and selections of ends, with applications to grid graphs and to measurable and Baire measurable edge weightings. Finally, Section~\ref{sec:questions} presents further questions.
\section{Preliminaries}\label{sec:preliminaries}
\subsection{Basics of graph theory and descriptive set theory}

Let $G=(X,E)$ be an undirected simple graph, where $X$ is a set and $E$ is a symmetric and irreflexive subset of $X\times X$. The sets $X$ and $E$ are called the \textit{vertex set} and the \textit{edge set}, respectively. For convenience, we sometimes write $xy$ instead of $(x,y)\in E$. We say that $x$ and $y$ are \textit{incident} if $xy\in E$, and an edge $xy$ is \textit{incident} to its endpoints $x$ and $y$.  A \textit{cycle} of length $n\geq3$ consists of distinct vertices $x_0,\ldots,x_{n-1}$ such that $x_ix_{i+1}\in E$ for $0\leq i<n-1$ and $x_{n-1}x_0\in E$. A graph is \textit{acyclic}, or a \textit{forest}, if it contains no cycles. The graph $G$ is \textit{bipartite} if its vertex set has a partition $X=A\sqcup B$ such that every edge connecting one vertex in $A$ and one vertex in $B$, equivalently, $G$ contains no odd cycle.

A \textit{path} of length $n$ is a sequence of distinct vertices $x_0,\ldots,x_n$ such that $x_ix_{i+1}\in E$ for every $i<n$. Two vertices belong to the same \textit{connected component} if they can be joined by a path. A connected acyclic graph is called a \textit{tree}. A \textit{spanning subgraph} of $G$ is a graph $H=(X,E_H)$ with $E_H\subseteq E$. It is a \textit{spanning subforest} if it is acyclic. If $G$ is connected, a connected spanning subforest is called a \textit{spanning tree}.

Let $M\subseteq X$. The \textit{induced subgraph} on $M$, denoted by $G[M]$, has vertex set $M$ and edge set $E\cap(M\times M)$. We say that $M$ is an \textit{independent set} if no two vertices in $M$ are incident. An independent set $M$ is \textit{maximal} if every vertex in $X\setminus M$ is incident to some vertex in $M$. Two disjoint subsets of $X$ are said to be \textit{incident} if some edge has one endpoint in each of them.

For each $x\in X$, the \textit{neighborhood} of $x$ is $N_G(x)=\{y\in X:xy\in E\}$. The \textit{degree} of $x$ is $\deg_G(x)=|N_G(x)|$. A graph is \textit{locally finite} if every vertex has finite degree, and \textit{locally countable} if every vertex has at most countably many neighbors. It is \textit{$d$-regular} if $\deg_G(x)=d$ for every $x\in X$. We write $\Delta(G)=\sup\{\deg_G(x):x\in X\}$. A vertex of degree $0$ is an \textit{isolated vertex}, and a vertex of degree $1$ is a \textit{leaf}. An \textit{isolated edge} is a connected component consisting of two incident vertices.

The \textit{path distance} $d_G(x,y)$ is the length of a shortest path from $x$ to $y$, with $d_G(x,y)=\infty$ if they lie in different components. For $A,B\subseteq X$, put $d_G(x,A)=\inf\{d_G(x,a):a\in A\}$ and $d_G(A,B)=\inf\{d_G(a,b):a\in A,\ b\in B\}$, where the infimum of the empty set is $\infty$. For $r\in\mathbb{N}$, write $B_G(A,r)=\{x\in X:d_G(x,A)\leq r\}$ and $B_G(x,r)=B_G(\{x\},r)$.

We also introduce some basic concepts in descriptive set theory and the theory of countable Borel equivalence relations; our references are~\cite{Kec95,Slu15,KM04}. A topological space is \textit{Polish} if it is separable and completely metrizable. A \textit{standard Borel space} is a measurable space $(X,\mathcal{B})$ such that $\mathcal{B}$ is the Borel $\sigma$-algebra of some Polish topology on $X$. Such a topology is called \textit{compatible} with the Borel structure. A map $f:X\longrightarrow Y$ between standard Borel spaces is called \textit{Borel} if $f^{-1}(A)$ is Borel for every Borel set $A\subseteq Y$. Borel subsets, countable products, and countable disjoint unions of standard Borel spaces are again standard Borel spaces, equipped with their natural Borel structures.

A map $f:X\longrightarrow Y$ is \textit{finite-to-one} if $f^{-1}(y)$ is finite for every $y\in Y$, and \textit{bounded-to-one} if there is an integer $m$ such that $|f^{-1}(y)|\leq m$ for every $y\in Y$. A map $f:X\longrightarrow X$ is \textit{aperiodic} if $f^n(x)\neq x$ for every $x\in X$ and every $n\geq1$.

An equivalence relation $E$ on a standard Borel space $X$ is a \textit{countable}, respectively \textit{finite}, \textit{Borel equivalence relation} if $E\subseteq X\times X$ is Borel and every equivalence class is countable, respectively finite. We write $[x]_E=\{y\in X:xEy\}$. A countable Borel equivalence relation is \textit{aperiodic} if all its classes are infinite. For $A\subseteq X$, its \textit{$E$-saturation} is $[A]_E=\{x\in X:\exists a\in A,\ xEa\}$. If $A$ is Borel, then $[A]_E$ is Borel by~\cite[Corollary~1.2.4]{Slu15}. A set $A$ is \textit{$E$-invariant} if $[A]_E=A$. A countable Borel equivalence relation $E$ is \textit{smooth} if it admits a \textit{Borel transversal}, that is, a Borel subset of $X$ meeting each equivalence class in exactly one point. Every finite Borel equivalence relation is smooth by~\cite[Proposition~1.4.4]{Slu15}. Let $\pi_E:X\longrightarrow X/E$ be the quotient map, and equip $X/E$ with the quotient $\sigma$-algebra consisting of those sets $A\subseteq X/E$ for which $\pi_E^{-1}(A)$ is Borel. By~\cite[Proposition~6.3]{KM04}, this quotient is standard Borel if and only if $E$ is smooth. In particular, the quotient by a finite Borel equivalence relation is standard Borel. A countable Borel equivalence relation $E$ on $X$ is \textit{hyperfinite} if there are finite Borel equivalence relations $F_n$ on $X$ such that $F_n\subseteq F_{n+1}$ for every $n\in\mathbb{N}$ and $E=\bigcup_{n\in\mathbb{N}}F_n$; see~\cite[Definition~3.1.1]{Slu15}.

Let $X$ be a standard Borel space. For each $n\geq0$, write $[X]^n=\{A\subseteq X:|A|=n\}$, and let $[X]^{<\infty}=\bigcup_{n\geq0}[X]^n$ be the space of all finite subsets of $X$. Fix a Borel linear order $<$ on $X$ and identify each $n$-element subset with its unique increasing enumeration. This identifies $[X]^n$ with the Borel subset $\{(x_1,\ldots,x_n)\in X^n:x_1<\cdots<x_n\}$ of $X^n$, so it is standard Borel. We equip $[X]^{<\infty}$ with the countable disjoint union Borel structure, as we do for the space $X^{<\infty}=\bigsqcup_{n\geq0}X^n$ of finite sequences.

We recall the Lusin--Novikov theorem, which will be used repeatedly to verify Borelness.

\begin{theorem}[Lusin--Novikov~{\cite[Theorem~18.10]{Kec95}}]
Let $X,Y$ be standard Borel spaces and let $E\subseteq X\times Y$ be Borel. If every section $E_x:=\{y\in Y:(x,y)\in E\}$ is countable, then there are Borel partial functions $f_n:D_n\longrightarrow Y$, with Borel domains $D_n\subseteq X$, such that
\[
E_x=\{f_n(x):n\in\mathbb{N},\ x\in D_n\}\qquad\text{for every }x\in X.
\]
\end{theorem}

Consequently, the projection of $E$ onto $X$ is Borel, and an element of each nonempty section can be chosen in a Borel way. For finite nonempty sections, one may take the least element with respect to a fixed Borel linear order on $Y$. Thus choices among nonempty finite sets of local configurations specified by Borel conditions can be made Borelly. We will use this theorem throughout the paper and omit routine Borelness verifications.

\subsection{Borel combinatorics}
We introduce some results in Borel combinatorics. A \emph{Borel graph} $G=(X,E)$ is a graph whose vertex set $X$ is a standard Borel space and whose edge set $E$ is a Borel subset of $X \times X$ which is equipped with the product $\sigma$-algebra. By~\cite[Corollary~13.4]{Kec95}, any Borel subset of a standard Borel space is again a standard Borel space with respect to the induced $\sigma$-algebra, so $E$ is a standard Borel space. It is also easy to see that if $D\subseteq X$ is a Borel subset, then the induced subgraph $G[D]$ is a Borel graph with the induced Borel structure.

The \textit{connectedness relation} of $G$, denoted by $E_G$, is defined by
\[
xE_Gy \iff x\text{ and }y\text{ belong to the same connected component}.
\]
For a locally finite Borel graph, $E_G$ is a countable Borel equivalence relation~\cite[Corollary~5.3]{Pik21}. We call $G$ smooth if $E_G$ is smooth. A set $A\subseteq X$ is $G$-invariant if it is $E_G$-invariant, equivalently, if it is a union of connected components of $G$.

Let $G=(X,E)$ be a Borel graph and $k$ be an integer. We let $d_G:X\times X \longrightarrow \mathbb{N}\cup \{\infty\}$ denote the distance function of the path metric of $G$. The $k$-\textit{power graph} of $G=(X,E)$, which we denote by $G^k$, is a graph with vertex set $X$ and edge set $\{(x,y):0<d_G(x,y)\le k\}$. By~\cite[Corollary~5.2]{Pik21}, if $G$ is a locally finite Borel graph, so is $G^k$.

We say that $w:E\longrightarrow\{1,\ldots,k\}$ is a \textit{Borel vertex-coloring k-edge weighting} if $w$ is a vertex-coloring k-edge weighting and is Borel measurable. Here we equip $k$ with the discrete $\sigma$-algebra. Suppose $Y$ is a standard Borel space. A \emph{Borel proper $Y$-coloring} of $G$ is a Borel map $c: X \longrightarrow Y$ such that for each $x,y\in X$, $xEy \implies c(x)\ne c(y)$. The Borel chromatic number of $G$ is defined by
$$ \chi_B(G)=\min\left\{|Y| \mid \text{There exists a Borel proper coloring } c: X \rightarrow Y\right\} $$
Every locally finite Borel graph $G$ satisfies $\chi_B(G)\leq\aleph_0$, and $\chi_B(G)\leq d+1$ if $\Delta(G)\leq d<\infty$~\cite[Propositions~4.5 and~4.6]{KST99}. For a locally countable Borel graph with countable Borel chromatic number, every Borel independent set extends to a Borel maximal independent set~\cite[Proposition~4.2]{KST99}. In particular, every locally finite Borel graph admits a Borel maximal independent set.

For a locally finite Borel graph $G$, let $[G]^{<\infty}=\{A\in[X]^{<\infty}:A\neq\varnothing\text{ and }G[A]\text{ is connected}\}$, a Borel subspace of $[X]^{<\infty}$. Given a Borel family $\mathcal{A}\subseteq[G]^{<\infty}$, its \textit{intersection graph} $I_{\mathcal{A}}$ has vertex set $\mathcal{A}$, with distinct $A,B\in\mathcal{A}$ adjacent exactly when $A\cap B\neq\varnothing$. This is a locally countable Borel graph, and $\chi_B(I_{\mathcal{A}})\leq\aleph_0$ by~\cite[Proposition~4.10]{KM20}, since each member of $\mathcal{A}$ is contained in an $E_G$-class. Consequently, $\mathcal{A}$ admits a maximal Borel subfamily of pairwise disjoint sets.

\subsection{Ends of Borel graphs}

Ends of Borel graphs and their connections with amenability, hyperfiniteness, and Borel selection problems have been studied in~\cite{Ada90,Mil09,HM09}. One-ended spanning forests and their applications to measurable combinatorics have been investigated in~\cite{CMT16,CGMT21,BPZ24}. We recall the relevant definitions below, following~\cite{Mil09,HM09,CMT16}.

Let $G=(X,E)$ be a locally finite Borel graph. A \textit{ray} of $G$ is an injective sequence $\alpha\in X^{\mathbb{N}}$ such that $(\alpha(n),\alpha(n+1))\in E$ for every $n\in\mathbb{N}$. The set of all rays is denoted by $[G]^\infty$. It is a Borel subset of $X^{\mathbb{N}}$, where $X^{\mathbb{N}}$ carries the product Borel structure, and hence inherits a standard Borel structure. Two rays $\alpha,\beta\in[G]^\infty$ are called \textit{end equivalent} if, for every finite set $S\subseteq X$, they can be connected by a path contained in $X\setminus S$. We denote this relation by $E_T$. An \textit{end} of $G$ is an equivalence class of $E_T$. We use the notation $E_T$ because if $G$ is acyclic, then end equivalence is exactly tail equivalence, which means that two rays eventually agree except on finite initial segments. For a cardinal $\kappa$, we say $G$ is \textit{$\kappa$-ended} if each component of $G$ has exactly $\kappa$ ends. When $\kappa$ is a finite integer, say $k$, being $k$-ended is equivalent to the following: for every component $C$ of $G$, removing any finite subset $S\subseteq X$ divides that component into at most $k$ infinite parts, and there exists a finite set $F\subseteq{X}$ whose removal divides that component into exactly $k$ infinite parts. Precisely, the induced subgraph $G[C\setminus S]$ has at most $k$ infinite components, and there exists a finite set $F\subseteq{X}$ such that $G[C\setminus F]$ has exactly $k$ infinite components.

A \textit{line} of $G$ is a set of distinct vertices $\{x_n:n\in\mathbb{Z}\}$ satisfying $(x_m,x_n)\in E\iff |m-n|=1$ for all $m,n\in\mathbb{Z}$. In an acyclic graph, two distinct ends in the same connected component determine a unique line. To see this, choose rays representing the two ends and remove finite initial segments so that the remaining rays are disjoint, and we connect them from shortest path.

Let $k$ be an positive integer. A \textit{Borel selection of $k$ ends from each component} is an $E_T$-invariant Borel set $A\subseteq[G]^\infty$ such that $A\cap[G[C]]^\infty$ is the union of exactly $k$-many $E_T$-equivalence classes for every connected component $C$ of $G$. More generally,  a \textit{Borel selection of a nonempty finite set of ends from each component} is an $E_T$-invariant Borel set $B\subseteq[G]^\infty$ such that $B\cap[G[C]]^\infty$ is the union of finitely (non-zero) many $E_T$-equivalence classes for every connected component $C$ of $G$.

A \textit{one-ended spanning forest} of $G$ is a forest $F=(X,E_F)$ with $E_F\subseteq E$ whose connected components each have exactly one end. It is Borel if $E_F$ is Borel. We say a partial function $p:D\longrightarrow X$ is \textit{contained in $G$} if $(x,p(x))\in E$ for every $x\in D$. We say that $p$ is \textit{one-ended} if there is no sequence $(x_n)_{n\in\mathbb{N}}$ in $X$ satisfying $x_{n+1}\in D$ and $p(x_{n+1})=x_n$ for every $n\in\mathbb{N}$. If $D$ and $p$ are Borel, we call $p$ a \textit{one-ended Borel partial function}. In particular, when $D=X$, a one-ended Borel function contained in $G$ generates a Borel one-ended spanning forest. Conversely, a Borel one-ended spanning forest determines a one-ended Borel function by sending each vertex to the next vertex on the unique ray starting there.

We will repeatedly use the following result of Conley, Marks, and Tucker-Drob to construct one-ended Borel partial functions.

\begin{proposition}[{\cite[Proposition~3.1]{CMT16}}]
Let $G=(X,E)$ be a locally finite Borel graph on a standard Borel space $X$, and let $A\subseteq X$ be Borel. There exists a one-ended Borel partial function $p:[A]_{E_G}\setminus A\longrightarrow [A]_{E_G}$ contained in $G$.
\end{proposition}

\section{Some general upper bounds}\label{sec:upper-bounds}

In this section, we will give some general bounds for some special classes of Borel graphs.

Note that verifying whether an edge weighting is vertex-coloring is a local problem with a finite color set. Indeed, we only need to check the ball of radius $2$. Thus all smooth graphs behave in the same way as finite graphs. We present some details as follows.

Precisely, a \textit{structured graph} is a pair $(G,\tau)$, where $\tau$ is a partial function from $X^{<\infty}$ to $\mathbb{N}$ such that for some $l\in \mathbb{N}$, every  $x\in dom(\sigma)$ is of length
at most $l$	. A \textit{local coloring problem} is a pair $\Pi=(r, P)$, where $r\in{\mathbb{N}}$ and $P$ is a map that sends isomorphism classes of finite labeled graphs to $\{0,1\}$. A $\Pi$-coloring of a locally finite structured graph $(G,\tau)$ is a function $\varphi : X\longrightarrow \mathbb{N}$ such that $P([G,\tau,\varphi ,x]_r)=1$ for all $x\in{X}$. For further details on these definitions, see~\cite[\S~2]{Ber23}.

\begin{proposition}\label{smooth}
Suppose $G$ is a smooth locally finite Borel graph with no isolated edges. Then $G$ admits a Borel vertex-coloring $3$-edge weighting.
\end{proposition}

\begin{proof}
Let $\mathcal{E}:=\{\{x,y\}:xy\in E(G)\}$, equipped with the Borel structure inherited from $[X]^2$. Consider the following incidence graph $I$, where the vertex set of $I$ is the disjoint union $X\sqcup\mathcal{E}$, and $x$ is incident to $e$ if and only if $x\in{e}$, for $x\in{X}$ and $e\in{\mathcal{E}}$. The graph $I$ is still a locally finite Borel graph. If $T$ is a Borel transversal for $G$, then it is easy to see that $T$ is also a Borel transversal for $I$, because each component of $I$ is obtained by adding some two-point sets to a component of $G$, and $T$ does not intersect these two-point sets. Thus $I$ is a locally finite smooth Borel graph.

Let $\tau:X\sqcup\mathcal{E} \longrightarrow \{\mathsf{v},\mathsf{e}\}$, where $\tau(x)=\mathsf{v}$ if $x\in{X}$ and $\tau(x)=\mathsf{e}$ if $x\in{\mathcal{E}}$. Then $(I,\tau)$ becomes a structured Borel graph. We define a local coloring problem $\Pi=(2,P)$. For any $\varphi:X\sqcup\mathcal{E}\longrightarrow\{1,2,3\}$ and any $x\in X\sqcup\mathcal{E}$, let $[I,\tau,\varphi,x]_2$ denote the isomorphism class of the rooted radius-$2$ ball around $x$, equipped with the restrictions of $\tau$ and $\varphi$. This ball is finite since $I$ is locally finite. Define $P([I,\tau,\varphi,x]_2)=1$ if and only if one of the following conditions holds:
\begin{enumerate}
\item $\tau(x)=\mathsf{v}$ and $\varphi(x)=1$;
\item $\tau(x)=\mathsf{e}$, $\varphi(x)\in\{1,2,3\}$, and $N_I(x)=\{u,v\}$ consists of two distinct vertices of type $\mathsf{v}$, all vertices in $N_I(u)\cup N_I(v)$ have type $\mathsf{e}$, and
\[
    \sum_{y\in N_I(u)}\varphi(y)
    \neq
    \sum_{y\in N_I(v)}\varphi(y).
    \]
\end{enumerate}
In all other cases, set $P([I,\tau,\varphi,x]_2)=0$. We also extend $P$ to all remaining isomorphism classes of finite rooted colored structured graphs by setting it equal to $0$. It is easy to see that a $\Pi$-coloring of $I$ is equivalent to a vertex-coloring $3$-edge weighting of $G$. So by~\cite[Theorem~2]{Sta23}, $G$ admits a vertex-coloring $3$-edge weighting. Then since $I$ is smooth and locally finite,~\cite[Theorem~1.2]{Ber22} gives a Borel vertex-coloring $3$-edge weighting
\end{proof}

\begin{remark}
So by the proceeding proposition, unless otherwise stated, we assume all the graph in this paper has only infinite connected components, since finite components of a graph is smooth.
\end{remark}

For locally finite Borel graphs, we show a general upper bound by directly adapting the algorithm of Kalkowski, Karo\'{n}ski, and Pfender~\cite{KKP10} to the Borel setting. 

\begin{theorem}
Suppose $G=(X, E)$ is a locally finite Borel graph. Then $G$ admits a Borel vertex-coloring $5$-edge weighting.
\end{theorem}
\begin{proof}
Let $D=\{x\in X: \deg(x)\ge 2\}$. Because there are no isolated edges, $D$ is nonempty. We apply~\cite[Proposition~4.2]{KST99} to the induced subgraph $G^2[D]$ to obtain a Borel maximal independent set $M$ such that $d_G(x,M)\le 2$ for each $x\in D$ and $d_G(a,b)> 2$ for distinct $a,b\in M$. Let $L_i=\{x\in X:d_G(x,M)=i\}, i=0,1,2,3$. $\{L_i:i=0,1,2,3\}$ indeed forms a Borel partition of $X$. To see this, for each $x\in X$, if $\deg(x)=1$, then there is a $y\in D$ such that $xEy$; otherwise, $xy$ would be an isolated edge. So $d_G(x,M)$ is at most $3$.

Our first goal is as follows: we will construct a Borel map $w:E\longrightarrow \{1,2,3,4,5\}$, a Borel map $h: X\setminus M\longrightarrow \mathbb{N}$ and a set $W(x)=\{h(x),h(x)+2\}$ such that the following conditions hold:
\begin{enumerate}
\item $h(x)\equiv 0,1 \pmod{4}$.
\item For every $x,y\notin M,$ if $xEy$, then $W(x)\cap W(y)=\emptyset$.
\item $s_w(x)=\sum_{zEx}w(zx)\in W(x)$.
\item $s_w(a)\ne s_w(z)$ for every $z\in N_G(a)$ and every $a\in M$.
\end{enumerate}

Now we begin the construction. By~\cite[Proposition~4.5]{KST99}, let $c:X\longrightarrow \mathbb{N}$ be a Borel proper coloring of $G^2$. For each $x\in L_i$, $i=1,2,3$, we choose a unique $p(x)$ as follows: if the set $Y=\{ y\in N_G(x)\cap L_i: c(y)> c(x)\}$ is nonempty, let $p(x)$ be the element of $Y$ with the least $c$-color. This element is unique because any two distinct vertices $u,v\in{N_G(x)}$ have distance at most $2$ and hence have distinct colors. Otherwise, let $p(x)$ be the unique element of $N_G(x)\cap L_{i-1}$ with the least $c$-color. Note that $N_G(x)\cap L_{i-1}$ is always nonempty because there is a path toward $L_0$, so $p$ is well-defined. Since $G$ is locally finite and $c$ is a Borel coloring, $p:X\setminus M\longrightarrow X$ is Borel.

Let $A_{i,j}=L_i\cap c^{-1}(j)$, where $i=1,2,3$ and $j\in \mathbb{N}$. $\{A_{i,j}\}_{j}$ forms a partition of $L_i$. We begin with $L_3$. By definition, each element of $L_3$ has degree $1$, and there is a unique edge connecting $L_3$ and $L_2$. We define $w_3(e)=3$ for every $e\in E$, $h_3(x)=1$ and $W_3(x)=\{1,3\}$ for each $x\in L_3$.

Then we turn to $L_2=\bigcup_{j}A_{2,j}$. For each $j\in \mathbb{N}$, let $P_2^j=L_3\cup \bigcup_{m<j}A_{2,m}$ and note that $P_2^0=L_3$. We will construct a tuple $(w_2^j,h_2^j,W_2^j)$ by induction, where $w_2^j: E\longrightarrow \{1,2,3,4,5\}$, $h_2^j:P_2^j\longrightarrow \mathbb{N}$ and $W_2^j(x)=\{h_2^j(x),h_2^j(x)+2\}$ with $h_2^j(x)\equiv 0,1 \pmod{4}$. When $j=0$, we have $P_2^j=L_3$, and we let $(w_2^0,h_2^0,W_2^0)$ be $(w_3,h_3,W_3)$.

Suppose we have defined $(w_2^j,h_2^j,W_2^j)$. Now, we define a Borel map $w_2^{j+1}:E\longrightarrow \{1,2,3,4,5\}$, a Borel map $h_2^{j+1}: P_2^{j+1}\longrightarrow \mathbb{N}$ and a set $W_2^{j+1}(x)=\{h_2^{j+1}(x),h_2^{j+1}(x)+2\}$ such that the following requirements are satisfied:
\begin{enumerate}
\item $w_2^{j+1}(e)=w_2^{j}(e)$ for each $e$ not incident to $A_{2,j}$, $h_2^{j+1}$ and $h_2^{j}$ are compatible \\ and $W_2^{j+1}(x)=W_2^{j}(x)$ for all $x\in P_2^{j}$.
\item $h_2^{j+1}(x)\equiv 0,1 \pmod{4}$.
\item For every $x,y\in P_2^{j+1},$ if $xEy$, then $W_2^{j+1}(x)\cap W_2^{j+1}(y)=\emptyset$.
\item $s_{w_2^{j+1}}(x)=\sum_{zEx}w_2^{j+1}(zx)\in W_2^{j+1}(x)$ for every $x\in P_2^{j+1}$.
\item For each $x\in P_2^{j+1}, y\notin P_2^{j+1}$, if $xEy$ and $y\ne p(x)$, then $w_2^{j+1}(xy)=3$.
\item For each $x\in P_2^{j+1}$ such that $p(x)\notin P_2^{j+1}$:
\begin{enumerate}
\item If $s_{w_2^{j+1}}(x)=h_2^{j+1}(x),$ then $w_2^{j+1}(xp(x))\in\{2,3\}$.
\item If $s_{w_2^{j+1}}(x)=h_2^{j+1}(x)+2,$ then $w_2^{j+1}(xp(x))\in\{3,4\}$.
\end{enumerate}
\end{enumerate}
Since $c$ is a proper coloring of $G^2$, our construction does not change any edge or vertex twice at any step and is therefore well-defined.

To see why this construction can be done, we fix $x\in A_{2,j}$ and let $B_j(x)=N_G(x)\cap P_2^j$, $d_j(x)=|B_j(x)|$. Each $y\in B_j(x)$ satisfies our induction conditions. At step $j+1$, we will change the value of $w_2^j(yx)$, keep the value of $w_2^j(zy)$ unchanged for each $z\ne x, zEy$, and ensure that $s_{w_2^{j+1}}(y)\in W_2^j(y)$. By our induction conditions, there are the following four cases for the value of $(s_{w_2^j}(y),w_2^{j}(yx))$, with the following choices:
\begin{enumerate}
\item If $s_{w_2^j}(y)=h_2^j(y)$ and $w_2^{j}(yx)=3$, then we have two possible choices for $w_2^{j+1}(yx)$: we can still keep $w_2^{j+1}(yx)=3$ or let $w_2^{j+1}(yx)=5$, which ensures that $s_{w_2^{j+1}}(y)\in W_2^j(y)$.
\item If $s_{w_2^j}(y)=h_2^j(y)+2$ and $w_2^{j}(yx)=3$, then we can let $w_2^{j+1}(yx)=1$ or $w_2^{j+1}(yx)=3$.
\item If $s_{w_2^j}(y)=h_2^j(y)$ and $w_2^{j}(yx)=2$, then we can let $w_2^{j+1}(yx)=2$ or $w_2^{j+1}(yx)=4$.
\item If $s_{w_2^j}(y)=h_2^j(y)+2$ and $w_2^{j}(yx)=4$, then we can let $w_2^{j+1}(yx)=2$ or $w_2^{j+1}(yx)=4$.
\end{enumerate}
We give some explanation here. Consider the edge $yx$. If $x\ne p(y)$, then, by the induction conditions, we have $w_2^j(yx)=3$, which leads to case $1$ or $2$. If $x=p(y)$, then $w_2^j(yx)=2,3,4$ (by induction conditions $6.(a)(b)).$ The case $w_2^j(yx)=3$ again leads to case $1$ or $2$, and $w_2^j(yx)=2,4$ leads to case $3$ or $4$. No other case is possible.

For each $y$, we denote the smaller possible value of $w_2^{j+1}(yx) $ by $l_y$; the larger one is $l_y+2$. Now for each edge $yx$, where $y\in B_j(x)$, we first let $w_2^{j+1}(yx)=l_y$ and let $w_2^{j+1}(xp(x))=2$, and the other edges $zx$ for $z\in N_G(x)\setminus (B_j(x)\cup \{p(x)\})$ retain the value $3$. In this case, we have
\begin{center}
$s_{w_2^{j+1}}(x)=\sum_{y\in B_j(x)}l_y+2+3\bigl(\deg(x)-d_j(x)-1).$
\end{center}

This is the least possible value for $s_{w_2^{j+1}}(x)$ and we denote it by $a_j(x)$. There are two valid values of $w_2^{j+1}(yx)$ for each $y\in B_j(x)$, that is, $l_y$ and $l_y+2$. Moreover, we can choose the value of $w_2^{j+1}(xp(x))$ from $\{2,3,4\}$. All these choices can be made in a Borel way, and all possible values of $s_{w_2^{j+1}}(x)$ form an integer interval $[a_j(x),a_j(x)+2d_j(x)+2]_\mathbb{Z}$.

Fix one possible value of $s_{w_2^{j+1}}(x)$ and denote it by $t$. We claim that $t$ uniquely determines $h_2^{j+1}(x)$ and hence $W_2^{j+1}(x)$. To ensure that induction condition 2 holds, if $t\equiv 2,3 \pmod{4}$, we let $h_2^{j+1}(x)=t-2$. If $t\equiv 0,1 \pmod{4}$, we let $h_2^{j+1}(x)=t$. We denote this unique value by $b(t)$. Then we will choose a $t$ such that $\{b(t),b(t)+2\}\cap W_2^j(y)=\emptyset$ for every $y\in B_j(x)$. To see this, note that $\{b(t),b(t)+2\}\cap W_2^j(y)\ne\emptyset \iff t=h_2^{j}(y)$ or $t=h_2^{j}(y)+2$. So each $y$ will block at most two possible values of $t$, and at most $2d_j(x)$ values in $[a_j(x),a_j(x)+2d_j(x)+2]_\mathbb{Z}$ will be blocked. We also note that if $t \in[a_j(x)+1,a_j(x)+2d_j(x)+1]_\mathbb{Z}$, then it will satisfy induction condition 6. In fact, if we choose $t=a_j(x)$, this means $w_2^{j+1}(xp(x))$ can only be $2$, which contradicts $6.(b)$ if $h_2^{j+1}(x)=s_{w_2^{j+1}}(x)-2$. For a similar reason, $a_j(x)+2d_j(x)+2$ might be an invalid value of $t$. It is easy to see that every other choice is valid: if $t-a_j(x)$ is odd, let $w_2^{j+1}(xp(x))=3$. If $t-a_j(x)$ is even and $t=h_2^{j+1}+2$, let $w_2^{j+1}(xp(x))=4$. If $t-a_j(x)$ is even and $t=h_2^{j+1}$, let $w_2^{j+1}(xp(x))=2$.

Now at most $2d_j(x)+2$ points, including the two endpoints, might be blocked. However, we have $2d_j(x)+3$ choices, so we can always choose the least value that is not blocked from $[a_j(x)+1,a_j(x)+2d_j(x)+1]_\mathbb{Z}$ for $s_{w_2^{j+1}}(x)$. This value satisfies all induction conditions. Thus the induction is complete, and we obtain a sequence $(w_2^j,h_2^j,W_2^j)_{j\in\mathbb{N}}$. Now let $w_2=\lim_j w_2^j$, $h_2=\bigcup_j h_2^j$, $W_2(x)=\{h_2(x),h_2(x)+2\}$. The limit exists because these maps are compatible and each edge changes at most twice. Note that $w_2$ and $h_2$ are still Borel.

The construction of $(w_1,h_1,W_1)$ is the same; we only need to define the initial step $(w_1^0,h_1^0,W_1^0)$ as $(w_2,h_2,W_2)$. Let $(w,h,W)$=$(w_1,h_1,W_1)$. Then conditions $1,2,3$ of our goal are satisfied.

Now we deal with the Borel maximal independent set $M$. We continue to use the notation above. We will construct a Borel map $\widehat{w}: E\longrightarrow \{1,2,3,4,5\}$ such that $\widehat{w}$ is compatible with $w$ except on those edges incident to $M$. For these edges, we change their values to satisfy the edge-weighting requirement. This $\widehat{w}$ will be our final result. Note that there is no $p(x)$ for $x\in M$, and we no longer need the induction.

Fix $x\in{M}$, and choose the least possible value $w(yx)=l_y$ for each $y\in{N_G(x)}$. Let $k=\sum_{y\in{N_G(x)}}l_y$. We consider three cases to obtain $\widehat{w}$:
\begin{enumerate}
\item If $k\equiv 2,3 \pmod{4}$, then we let $\widehat{w}(yx)=l_y$ for all $y\in{N_G(x)}$.
\item If $k\equiv 0,1 \pmod{4}$ and there is a $y$ such that $h(y)\ne k$, then we choose a unique such $y$ in a Borel way, denote it by $y_0$, and then let $\widehat{w}(y_0x)=l_{y_0}+2$ for this $y_0$.
\item If $k\equiv 0,1 \pmod{4}$ and there is no $y$ such that $h(y)\ne k$, then we arbitrarily choose two vertices $y_0$, $y_1 \in{N_G(x)}$ in a Borel way, and then let $\widehat{w}(y_0x)=l_{y_0}+2$, $\widehat{w}(y_1x)=l_{y_1}+2$ for $y_0$, $y_1$.
\end{enumerate}
The operation in case $3$ can be carried out since each $x\in{M}$ has $\deg(x)\ge 2$. In case $1$, by the inductive construction on $L_1$, when $\widehat{w}(yx)$ takes the value $l_y$, we have $s_{\widehat{w}}(y)\equiv 0,1 \pmod{4}$, which is different from $k$. In case $2$, we increase $\widehat{w}(y_0x)$ and $k$ at the same time while keeping $s_{\widehat{w}}(y_0)\ne k+2$. Now $k+2\equiv 2,3 \pmod{4}$, which differs from $s_{\widehat{w}}(y)$ for every $y\in{N_G(x)}\setminus\{y_0\}$ for the same reason as in case $1$. Case $3$ is obvious. Note that $M$ is a $G^{2}$-independent set, so for each $x\in{X\setminus M}$, at most one edge incident to both this vertex and M changes. Thus these changes are well-defined. Finally, we obtain $\widehat{w}$, which is clearly a Borel vertex-coloring $5$-edge weighting.
\end{proof}

Adapting the algorithm in~\cite[Section~3]{Prz21} to the Borel setting, we also obtain the following upper bound for Borel regular graphs.

\begin{theorem}\label{regular}
Suppose $G=(X, E)$ is a regular locally finite Borel graph. Then $G$ admits a Borel vertex-coloring $4$-edge weighting.
\end{theorem}

\begin{proof}
First, we partition $X$ into three parts, as in the finite proof~\cite[Section~3]{Prz21}. Let $I$ be a Borel maximal independent set of $G$, and put
\[
 R_1:=\{x\in X\setminus I:N_G(x)\subseteq I\},
 \qquad
 R_2:=X\setminus(I\cup R_1).
\]
Every vertex outside $I$ has a neighbor in $I$. Moreover, $R_1$ is independent, there are no edges between $R_1$ and $R_2$, and the restriction $G[R_2]$ of $G$ to $R_2$ has no isolated vertices.

To run the algorithm Borelly, we partition $R_2$ into smaller pieces. Take a Borel maximal independent set $A_0$ of $(G[R_2])^2$, and let
\[
 A_1:=\{x\in R_2:d_{G[R_2]}(x,A_0)=1\},
 \qquad
 A_2:=\{x\in R_2:d_{G[R_2]}(x,A_0)=2\}.
\]
By maximality, $R_2=A_0\sqcup A_1\sqcup A_2$.

For $\ell=1,2$, let $G^2[A_\ell]$ denote the graph on $A_\ell$ in which distinct $x,y$ are incident when $d_{G[R_2]}(x,y)\leq2$. Thus the square is taken before restricting to the layer. Take finite Borel colorings $c_1,c_2$ of these graphs, with color classes
\[
 A_1=\bigsqcup_{i<N}A_1^i,
 \qquad
 A_2=\bigsqcup_{j<M}A_2^j.
\]
Distinct vertices in any one class have distance greater than $2$ in $G[R_2]$. The same holds for distinct vertices of $A_0$. All auxiliary graphs have bounded degree, so the independent sets and finite Borel colorings used here exist.

We process the classes in the order
\begin{equation}\label{eq:order}
 A_2^0\longrightarrow\cdots\longrightarrow A_2^{M-1}
 \longrightarrow A_1^0\longrightarrow\cdots\longrightarrow A_1^{N-1}
 \longrightarrow A_0,
\end{equation}
processing all vertices of each class simultaneously. Every vertex in $A_2$ has a neighbor in $A_1$, and every vertex in $A_1$ has a neighbor in $A_0$. Thus every vertex outside $A_0$ has a later neighbor.

For every $x\in A_1\cup A_2$, choose a supporting edge $e_x^I$ joining $x$ to $I$, and a first forward edge $e_x^{R_2}\in G[R_2]$. Fix a Borel linear order on $X$, and choose $e_x^I$ using the least neighbor of $x$ in $I$.

The first forward edge is chosen according to~\eqref{eq:order}. For $x\in A_2^{j_x}$, first look for a neighbor in a class $A_2^j$ with $j>j_x$, taking the least such $j$. If there is none, take a neighbor in $A_1^i$ for the least possible $i<N$. For $x\in A_1^{i_x}$, similarly take a neighbor in the least class $A_1^i$ with $i>i_x$, if possible; otherwise take its neighbor in $A_0$. No two neighbors of a vertex lie in the same class in~\eqref{eq:order}, so these forward neighbors are unique. In particular, the choices are Borel. Vertices of $A_0$ have no designated supporting or forward edges.

Define an initial weighting $w$ by
\[
w(e):=
\begin{cases}
1,&\text{if }e=e_x^{R_2}\text{ for some }x\in A_1\cup A_2,\\
4,&\text{if }e=e_x^I\text{ for some }x\in A_1\cup A_2,\\
2,&\text{if }e\in G[R_2]\text{ is not a first forward edge},\\
3,&\text{otherwise}.
\end{cases}
\]
We continue to write $w$ for the current weighting as it is modified. Throughout the procedure on $R_2$,
\[
\begin{aligned}
 w(e)&\in\{1,2,3\}&&\text{for }e\in G[R_2],\\
 w(e)&\in\{3,4\}&&\text{for edges joining }I\text{ to }R_2,\\
 w(e)&=3&&\text{for edges joining }I\text{ to }R_1.
\end{aligned}
\]
In the final procedure on $R_1$, the last range will become $\{2,3,4\}$. Our final goals are
\begin{equation}\label{eq:goals}
\begin{aligned}
 s_w(x)&\in\{3d-1,4d\}&&\text{for }x\in R_1,\\
 s_w(x)&\geq3d&&\text{for }x\in I,\\
 s_w(x)&<4d&&\text{for }x\in I\text{ with a neighbor in }R_1,\\
 s_w(x)&<3d&&\text{for }x\in R_2,
\end{aligned}
\end{equation}
together with distinct sums at incident vertices of $R_2$.

After a vertex is processed, its sum will remain fixed, although the weights of its incident edges may change later. For stage labels $(k,\ell)$, write $(k,\ell)<(k_0,\ell_0)$ when $A_\ell^k$ occurs before $A_{\ell_0}^{k_0}$ in~\eqref{eq:order}. Suppose that all earlier stages have been completed, and let $x\in A_{\ell_0}^{k_0}$. Define
\[
 N^-_{G[R_2]}(x):=
 \bigl\{y\in R_2:yG[R_2]x\text{ and }
       y\in A_\ell^k\text{ for some }(k,\ell)<(k_0,\ell_0)\bigr\}.
\]
If this set is empty, leave the weighting unchanged at $x$. Otherwise, use the following claim, which is the local adjustment from the finite proof.

\begin{claim}
By changing only the weights of $(y,x)$ and $e_y^I$, for $y\in N^-_{G[R_2]}(x)$, we can make $s_w(x)$ different from $s_w(y)$ for every such $y$, while preserving the sums at all previously processed vertices of $R_2$ and the prescribed edge-weight ranges. The choices can be made Borelly.
\end{claim}

\begin{proof}[Proof of the claim]
Let $y\in N^-_{G[R_2]}(x)$. The edge $(y,x)$ has not yet been modified: an edge of $G[R_2]$ is modified only when its later endpoint is processed.

If $(y,x)\neq e_y^{R_2}$, its current weight is $2$, and $w(e_y^I)\in\{3,4\}$. We may therefore make the appropriate one of the following changes:
\[
\begin{aligned}
 (w(y,x),w(e_y^I))&=(2,3)\longmapsto(1,4),\\
 (w(y,x),w(e_y^I))&=(2,4)\longmapsto(3,3).
\end{aligned}
\]
If $(y,x)=e_y^{R_2}$, then $x$ is the first later neighbor of $y$ to be processed. Thus neither this edge nor $e_y^I$ has yet been modified, and we may make the change
\[
 (w(y,x),w(e_y^I))=(1,4)\longmapsto(2,3).
\]
In all three cases we may also leave both weights unchanged. Each change preserves $s_w(y)$, stays within the required ranges, and changes $s_w(x)$ by $1$ or $-1$. No other previously processed vertex of $R_2$ is affected.

Put $b=|N^-_{G[R_2]}(x)|$. These $b$ optional changes concern distinct edges and distinct supporting edges, so they can be chosen independently. If $p$ of the available changes are increases and $q$ are decreases, where $p+q=b$, then the possible sums at $x$ include every integer from its current sum minus $q$ to its current sum plus $p$. Hence there are $b+1$ available values, while the $b$ earlier neighbors forbid at most $b$ values. At least one choice works. Since the graph has bounded degree, we can choose an admissible pattern in a Borel way.
\end{proof}

\smallskip The claim can be applied simultaneously at all vertices of the current class. Indeed, two such vertices are neither incident nor have a common neighbor in $G[R_2]$. Their modifications therefore use distinct internal edges and distinct supporting edges. Supporting edges from different operations may meet in $I$, but this causes no conflict: no sum in $I$ is fixed at this stage, and every incident edge retains weight at least $3$.

After all the classes in $A_2$ and $A_1$ have been processed, apply the same claim simultaneously to $x\in A_0$, now with $N^-_{G[R_2]}(x)=N_{G[R_2]}(x)$. This set is nonempty because $G[R_2]$ has no isolated vertices. Again choose the least admissible sum. All incident vertices of $R_2$ now have distinct sums.

We verify the bounds in~\eqref{eq:goals} for this stage. Every edge incident with $I$ has weight $3$ or $4$, so its vertices have sum at least $3d$. If a vertex of $I$ has a neighbor in $R_1$, the edge to that neighbor still has weight $3$, so its sum is less than $4d$.

For $x\in A_1\cup A_2$, immediately after processing $x$, its first forward edge still has weight $1$, its supporting edge still has weight $4$, and every other incident edge has weight at most $3$. Hence
\[
 s_w(x)\leq1+4+3(d-2)=3d-1.
\]
This sum is preserved by all later steps.

For $x\in A_0$, every incident edge has weight at most $3$, so $s_w(x)\leq3d$. Suppose that the least admissible choice at $x$ gave $s_w(x)=3d$. Then every edge incident with $x$ would have weight $3$. Consequently, for each $y\in N_{G[R_2]}(x)$, the edge $(y,x)$ is not $e_y^{R_2}$: a first forward edge starts at $1$ and can end at most at $2$. At the stage when $y$ was processed, the distinct edges $e_y^{R_2}$, $(y,x)$, and $e_y^I$ had weights $1$, $2$, and $4$, respectively, and every other incident edge had weight at most $3$. Thus its fixed sum satisfies
\[
 s_w(y)\leq1+2+4+3(d-3)=3d-2.
\]
Moreover, obtaining weight $3$ on $(y,x)$ required the change $(2,4)\mapsto(3,3)$ on $((y,x),e_y^I)$. Omitting this change for one neighbor $y$ would preserve every earlier sum and give $s_w(x)=3d-1$, distinct from all its neighbors' sums. This contradicts the choice of the least admissible sum. Therefore $s_w(x)<3d$ also for $x\in A_0$.

Finally, we deal with the edges between $I$ and $R_1$. Define the graph $G^2[R_1]$ on $R_1$ by
\[
 xG^2[R_1]y\quad\Longleftrightarrow\quad
 x\neq y\text{ and }d_G(x,y)=2.
\]
Take a finite Borel coloring $C_{R_1}$ of this graph, inducing a partition
\[
 R_1=\bigsqcup_{\ell<L}B_\ell.
\]
Process the classes in the order $B_0\longrightarrow B_1\longrightarrow\cdots\longrightarrow B_{L-1}$. For $x\in B_\ell$, all vertices in $N_G(x)$ lie in $I$. At the start of its stage, every edge incident with $x$ still has weight $3$.

If some $y\in N_G(x)$ has $s_w(y)\geq3d+1$, choose the least such $y$ and decrease $w(x,y)$ from $3$ to $2$. Then $s_w(x)=3d-1$, while $s_w(y)$ remains at least $3d$. Otherwise every $y\in N_G(x)$ has $s_w(y)=3d$. In this case, increase $w(x,y)$ from $3$ to $4$ for every $y\in N_G(x)$. Then $s_w(x)=4d$, and each neighbor's sum becomes $3d+1<4d$, since $d\geq2$.

Distinct vertices of $B_\ell$ have disjoint neighborhoods, so the operations within a class can be performed simultaneously. Throughout these stages,
\[
 3d\leq s_w(y)<4d
 \quad\text{for every }y\in I\text{ with a neighbor in }R_1.
\]
Thus the sums already fixed in $R_1$, each equal to $3d-1$ or $4d$, remain different from their neighbors' sums. No sum in $R_2$ is changed.

The final weighting satisfies~\eqref{eq:goals}. Edges within $R_2$ have already been distinguished; edges from $R_2$ to $I$ are distinguished by the threshold $3d$; and edges from $R_1$ to $I$ are distinguished by the preceding inequalities. There are no other types of edges. Finally, the entire construction uses finitely many Borel stages and Borel choices, so the resulting weighting $w$ is Borel.
\end{proof}

Recall that a triangle admits no vertex-coloring $2$-edge weighting. The following two examples show that, in the Borel setting, even acyclic graphs may fail to admit such a weighting.

\begin{example}
Let $F(2^{\mathbb Z})$ denote the free part of the Bernoulli shift action $\mathbb Z\curvearrowright 2^{\mathbb Z}$, and let $G$ be the associated Schreier graph with respect to the standard generators $\{-1,1\}$. Then $G$ is an acyclic $2$-regular Borel graph with $\chi_B(G)=3$~\cite[p.~3]{GJKS24}. Suppose that $w$ is a Borel vertex-coloring $2$-edge weighting of $G$. Since $G$ is $2$-regular, $s_w(x)\in\{2,3,4\}$ for every vertex $x$. Moreover, a vertex with sum $2$ cannot be adjacent to a vertex with sum $4$, since their common edge would have to receive both weights $1$ and $2$. As adjacent vertices have distinct sums, every edge has exactly one endpoint in the Borel set $A=\{x:s_w(x)=3\}$. Thus $A$ and its complement form a Borel bipartition of $G$, contradicting $\chi_B(G)=3$.
\end{example}

\begin{example}\label{F2}
A similar obstruction occurs for a $4$-regular acyclic Borel graph arising from the Bernoulli shift of the free group. Let $\mathbb F_2=\langle a,b\rangle$, and let $F(2^{\mathbb F_2})$ denote the free part of the Bernoulli shift action $\mathbb F_2\curvearrowright 2^{\mathbb F_2}$. Let $G$ be the associated Schreier graph with respect to the standard generators $\{a,a^{-1},b,b^{-1}\}$. Then $\chi_B(G)=5$~\cite[Theorem~1.2]{Mar16}. Suppose that $w$ is a Borel vertex-coloring $2$-edge weighting of $G$. Since $G$ is $4$-regular, $s_w(x)\in\{4,5,6,7,8\}$ for every vertex $x$. Moreover, a vertex with sum $4$ cannot be adjacent to a vertex with sum $8$, since their common edge would have to receive both weights $1$ and $2$. Thus identifying the sum values $4$ and $8$, while keeping $5$, $6$, and $7$ distinct, gives a Borel proper $4$-coloring of $G$, contradicting $\chi_B(G)=5$.
\end{example}

\section{Graphs with definable vertex-coloring 3-edge weightings}\label{sec:weightings}

In this section, we show that several basic examples admit Borel vertex-coloring $3$-edge weightings, including Schreier graphs induced by free Borel actions of $\mathbb{Z}^n$ and graphs generated by a single Borel function. Our main techniques involve using a fixed Borel vertex-$3$-coloring and analyzing the number of ends of the graphs. At the end of the section, we also discuss applications to $\mu$-measurable and Baire measurable edge weightings.

\subsection{\texorpdfstring{$k$}{k}-colorable graphs and one-ended spanning forests}
As recalled in Theorem~\ref{kcolor}, every finite $k$-colorable graph without isolated edges admits a vertex-coloring $k$-edge weighting when $k$ is odd. The main idea is to control the sums modulo a suitable integer. We begin with an elementary balance lemma.

\begin{lemma}\label{eq}
Let $G=(X=A\sqcup B,E)$ be a finite connected bipartite graph, and let $c:X\longrightarrow\{1,\ldots,k\}$ be a map. Then the following are equivalent:
\begin{enumerate}
    \item There exists an edge weighting $w:E\longrightarrow\{1,\ldots,k\}$ such that $s_w(x)\equiv c(x) \pmod{k}.$
    \item $\sum_{x\in A}c(x)\equiv \sum_{x\in{B}}c(x) \pmod{k}$.
\end{enumerate}
\end{lemma}
\begin{proof}
$1\implies 2$: Let $w$ be such an edge weighting. It is easy to see that $\sum_{x\in{B}}s_w(x)=\sum_{x\in A}s_w(x)$, because each edge weight contributes the same value to two vertices belonging to the two parts.

$2\implies 1$:  We list the elements of $X$ as $\{x_1,\dots,x_n\}$ and let $T$ be a least spanning tree of $G$ which rooted in $x_n$, such that $\forall i=1,\dots,n-1$, $\exists j>i$ such that $x_i$ is incident to $x_j$ in $T$. We first assign weight $k$ to every edge not appearing in $T$. We then define the values of $w$ on the other edges by induction.

We first arbitrarily assign weights to the edges incident to $x_1$ so that $s_w(x_1)\equiv c(x_1) \pmod{k}$. Suppose every edge incident to $x_{i-1}$ has been assigned a value of $w$. Now we consider the edges incident to $x_i$. Since there is an edge connecting it to some $x_j$ with $j>i$, we choose the least such $j$ and use the edge $x_ix_j$ to change the value of $s_w(x_i)$. Because $x_ix_j$ has not been assigned a value before, we can carry out this induction.

Finally, we deal with $x_n$. Without loss of generality, we may assume $x_n\in{A}$. By condition $2$, we have
$$c(x_n)\equiv \sum_{x\in{B}}c(x)-\sum_{x\in A\setminus \{x_n\}}c(x)\equiv \sum_{x\in{B}}s_w(x)-\sum_{x\in A\setminus \{x_n\}}s_w(x) \equiv s_w(x_n)\pmod{k}.$$
By induction, $w$ satisfies our requirement.
\end{proof}

This lemma transforms the weighting problem into a balance problem. We have the following theorem in Borel context.

\begin{lemma}\label{partial}
Let $G=(X,E)$ be a locally finite Borel graph, let $D\subseteq X$ be Borel, and let $p:D\longrightarrow X$ be a one-ended Borel partial function contained in $G$. For every integer $k\geq1$ and every Borel map $c:D\longrightarrow \{1,\ldots,k\}$, there is a Borel edge weighting $w:E\longrightarrow\{1,\ldots,k\}$ such that $s_w(x)\equiv c(x)\pmod{k}$ for every $x\in D$.
\end{lemma}

\begin{proof}
Let $p^{-1}(\{x\})=\{y\in D:p(y)=x\}$ for each $x\in D$. We construct Borel weightings $w_n:E\longrightarrow\{1,\ldots,k\}$ by induction. First, let $w_0(e)=k$ for every $e\in E$. Suppose $w_{n-1}$ has been defined. Let $w_n(e)=k$ for all outside edges. For each $x\in D$, choose the unique value $w_n(xp(x))\in\{1,\ldots,k\}$ satisfying
$$
w_n(xp(x))\equiv c(x)-\sum_{y\in p^{-1}(\{x\})}w_{n-1}(yx)\pmod{k}.
$$
We obtain a sequence of Borel edge weightings. We claim that the weight of each edge eventually stabilizes. Suppose that the weight of some edge $xp(x)$ changes infinitely many times. Then there must be some $y\in p^{-1}(\{x\})$ such that the weight of $yx$ also changes infinitely many times. Otherwise, by local finiteness, the sum over $p^{-1}(\{x\})$ would eventually stabilize, and so would the weight of $xp(x)$. Repeating this argument gives an infinite backward sequence, contradicting the one-endedness of $p$.

Finally, let $w=\lim_n w_n$. This limit exists, and $w$ is Borel. All outside incident edges have weight $k$, so $s_w(x)\equiv c(x)\pmod{k}$ for every $x\in D$, as required.
\end{proof}

Now we ready to prove the first main theorem in this section.

\begin{theorem}\label{oneend}
Suppose $G=(X, E)$ is a locally finite Borel graph. If $G$ is Borel $k$-colorable, then $G$ has a Borel vertex-coloring $(k+1)$-edge weighting. Moreover, if $G$ admits a Borel one-ended spanning forest, then $G$ has a Borel vertex-coloring $k$-edge weighting.
\end{theorem}

\begin{proof}
Suppose $c:X\longrightarrow \{1,\dots,k\}$ is a Borel proper coloring of $G$. We first extend the set $c^{-1}({k})$ to a Borel maximal independent set by $M$~\cite[Proposition~4.2]{KST99}. Then we change the colors of the vertices in $M\setminus c^{-1}({k})$ to $k$. Since no two vertices in $M$ are incident and no vertex outside $M$ has color $k$, the modified $c$ is still a Borel proper coloring of $G$. It is easy to see that $M$ is now a Borel maximal independent set of $G$ consisting of all vertices with color $k$.

We define a map $p:X\longrightarrow M$, where $p(x)$ is the vertex in $M$ nearest to $x$. If more than one vertex in $M$ is at the same minimal distance from $x$, we choose the least one with respect to the Borel linear order on $X$. Hence $p$ is a Borel map, and $\mathcal{P}=\{p^{-1}(x):x\in{M}\}$ gives a Borel partition of $X$. Let $E_T=\{(x,y)\in{E}:p(x)=y$ or $p(y)=x \}$. Then $(X,E_T)$ is a Borel graph in which each component is a finite tree rooted at an element of $M$.

We connect these small trees into larger trees in a Borel way so that each larger tree contains at least $2k$ small trees. Precisely, consider the Borel equivalence relation $F$ induced by this partition, that is, $x F y\iff p(x)=p(y)$. Because $M$ is a maximal independent set, $d_G(x,p(x))\le 1$. Together with the local finiteness of $G$, this shows that $F$ is a finite equivalence relation and hence smooth. Thus the quotient space $X/F$ is still standard Borel. This induces a Borel quotient graph $H=(X/F, \tilde{E})$, where $[x] \tilde{E} [y] \iff$ $[x]\ne [y]$ and $[x]$ and $[y]$ are incident in $G$. We again take a Borel maximal independent set $K$ of the $4k$-power graph of $H$ and similarly define $\tilde{p}:H\longrightarrow K$. Finally, for each $a\in{K}$, let $A=\bigcup \tilde{p}^{-1}(a)$, and note that $A$ is the union of the vertex sets of at least $2k$ small trees. Since these trees are connected in $G$, we choose the one with the least Borel color and restore exactly one incident edge from this tree to the other trees, step by step.

Now fix a larger tree. It is a finite bipartite graph; denote its two parts by $A$ and $B$. Without loss of generality, we assume that at least $k$ elements of $M$ lie in $A$. We can change the color of these elements from $k$ to $k+1$ to contribute $0,\dots,k$ to $\sum_{x\in A}c(x)$. Note that these changes preserve the fact that $c$ is a Borel proper coloring. Thus we can always ensure that the equation holds. By Lemma~\ref{eq}, this gives an edge weighting of the larger tree. We simply assign weight $k+1$ to every edge that is not in $E_T$ and has not been restored. We do this for every larger tree and obtain a Borel vertex-coloring $(k+1)$-edge weighting $w$.

For the moreover part, let $F$ be a Borel one-ended spanning forest of $G$, and let $q:X\longrightarrow X$ send each vertex to the next vertex on the unique ray of $F$ starting there. Then $q$ is a one-ended Borel function contained in $G$. Applying Lemma~\ref{partial} with $D=X$ and the original coloring $c$ (not change value), we obtain a Borel edge weighting $w:E\longrightarrow\{1,\dots,k\}$ such that $s_w(x)\equiv c(x)\pmod{k}$ for every $x\in X $. It is clear that $w$ is vertex-coloring.
\end{proof}

Note that $k$-colorability is not necessary. For example, let $[\mathbb{N}]^{\infty}$ be the standard Borel space of infinite subsets of $\mathbb{N}$, and let $f:[\mathbb{N}]^{\infty} \longrightarrow [\mathbb{N}]^{\infty}$ be the shift function defined by $f(A) = A \setminus \{\min(A)\}$. This graph is known as the \textit{shift graph}. By~\cite[Example~3.2]{KST99}, the Borel chromatic number of the shift graph is $\aleph_0$. However, this graph admits a Borel vertex-coloring $1$-edge weighting, because any two incident vertices have distinct degrees.

We conclude this subsection with an immediate consequence for Schreier graphs induced by free Borel $\mathbb{Z}^d$-actions.

\begin{corollary}\label{Zd}
Let $d\ge 2$, let $X$ be a standard Borel space, and let $\mathbb{Z}^{d}\curvearrowright X$ be a free Borel action, where $\mathbb{Z}^{d}$ is equipped with the standard generators $e_0,\dots,e_{d-1}$. Let $G$ be the Schreier graph induced by this action. Then $G$ admits a Borel vertex-coloring $3$-edge weighting.
\end{corollary}

\begin{proof}
By~\cite[Theorem~1.1]{GJKS24}, $G$ admits a Borel proper $3$-coloring. Moreover, for $d\ge 2$, since $\mathbb{Z}^{d}$ is one-ended and has polynomial growth,~\cite[Theorem~1.1(iii)]{BPZ24} implies that $G$ admits a Borel one-ended spanning forest (in fact, more strongly, admits one-ended tree).
\end{proof}

Note that the preceding corollary concerns $\mathbb{Z}^d$ for $d\geq2$, since the $\mathbb{Z}$-line is two-ended. The two-ended case is discussed in the next subsection.

\subsection{Two-ended bipartite graphs}

In this section we deal with the two-ended bipartite case. The proof idea are motivated by earlier works of ~\cite{Wei22} and ~\cite{Mil09}, which used to cut the whole graph into finite pieces in a Borel way. Combine with bipartiteness, we will construct a Borel 3-coloring such that our balance conditions hold.

\begin{proposition}\label{bipartite}
Suppose $G=(X,E)$ is a locally finite bipartite Borel graph. If $G$ is two-ended, then $G$ admits a Borel vertex-coloring $3$-edge weighting.
\end{proposition}
\begin{proof}
Let $\Phi\subseteq[G]^{<\infty}$ be the set of finite connected sets whose removal leaves exactly two infinite components in their connected component. As in~\cite{Wei22}, by modifying the proof of~\cite[Lemma~5.3]{Mil09}, there exists a maximal Borel set $\Psi\subseteq\Phi$ such that the path distance between any two distinct members of $\Psi$ is at least $10$. Fix such a $\Psi$.

Now, for each $M\in\Psi$, let $B(M)=\{x\in X:d_G(x,M)\leq1\}$ and consider the induced subgraph $G[B(M)]$. This is also a bipartite graph, and we write $B(M)=A_M\sqcup B_M$ for a bipartition. Note that this partition depends on $M$, and $G[B(M)]$ is finite, so this can be done in a Borel way. Let $\mathcal{A}=\{A_M:M\in\Psi\}$ and $W=\bigcup_{M\in\Psi}A_M$. We continue the definition in ~\cite{Wei22}, call distinct $S,T\in\Psi$ are $\mathcal{T}$-incident if there is a path from $S$ to $T$ which avoids all other points of $\bigcup\Psi$. Similarly, we use this terminology for distinct $S,T\in\mathcal{A}$ when there is a path from $S$ to $T$ which avoids all other points of $W$.

\begin{claim}\label{claim:2ended_1}
Each component of $G[X\setminus W]$ is incident to one or two elements of $\mathcal{A}$.
\end{claim}

\begin{proof}[Proof of claim \ref{claim:2ended_1}]
Each such component is incident to at least one element of $\mathcal{A}$, since the maximality of $\Psi$ implies that $W$ meets every component of $G$. Suppose that a component $U$ is incident to three distinct sets $A_{M_1},A_{M_2},A_{M_3}$. For every $M\in\Psi$, each vertex of $M\setminus A_M$ has all its neighbors in $A_M$ and is therefore isolated in $G[X\setminus W]$. Thus $U$ is disjoint from $\bigcup\Psi$. Since $A_{M_i}\subseteq B(M_i)$, each incidence of $U$ with $A_{M_i}$ extends to a path of length at most two from $U$ to $M_i$ which avoids $\bigcup\Psi$ except at its endpoint. Connecting these paths through $U$ shows that $M_1,M_2,M_3$ are pairwise $\mathcal{T}$-incident, contradicting ~\cite[Lemma~5.5]{Mil09}.
\end{proof}

\begin{claim}\label{claim:2ended_2}
$W$ is a Borel independent set, and each component of $G[X\setminus W]$ is finite.
\end{claim}

\begin{proof}[Proof of claim \ref{claim:2ended_2}]
Borelness is clear. To see that $W$ is independent, note that each $A_M$ consists of vertices from one bipartition class, and distinct $A_M,A_N$ have distance at least $8$. For finiteness, the preceding observation shows that every component of $G[X\setminus W]$ is either a singleton or contained in a component of $G[X\setminus\bigcup\Psi]$, which is finite by~\cite[Lemma~2.2]{Wei22}.
\end{proof}

Then we give $G$ a Borel coloring as follows: every vertex in $W$ receives color $3$. For every finite component of $G[X\setminus W]$, we similarly divide it into two parts and color them with $1$ and $2$. Denote this Borel coloring by $c$.

Let $\mathcal{S}_1$ denote the family of those finite components of $G[X\setminus W]$ which are incident to exactly one element of $\mathcal{A}$, and let $\mathcal{S}_2$ denote the others. We first deal with $\mathcal{S}_1$. Each $U\in\mathcal{S}_1$ is incident to a unique element $A\in\mathcal{A}$. Let $W_U=U\cup A$ and $D=\bigcup_{U\in\mathcal{S}_1}W_U$. Although different $W_U$ may overlap, each component of $G[D]$ is contained in the union of one $A_M$ and all members of $\mathcal{S}_1$ incident to it. This union is finite by local finiteness. Since every component of $G$ is infinite, each component of $G[D]$ has an edge to $X\setminus D$. Choose Borelly a spanning tree of each such finite component and an edge leaving it, and direct the tree toward the endpoint of this edge, which is sent outside $D$. This defines a one-ended Borel partial function $p:D\longrightarrow X$ contained in $G$. By Lemma~\ref{partial} and its construction, there exists a Borel edge weighting $w:E\longrightarrow\{1,2,3\}$ such that $s_w(x)\equiv c(x)\pmod{3}$ for every $x\in D$. We can assign weight $3$ to the edges outside $\{\{x,p(x)\}:x\in D\}$, and we will change their weights later.

Finally, we deal with $\mathcal{S}_2$. For each $U\in\mathcal{S}_2$, write $U=U_1\sqcup U_2$, where $U_1$ has been colored with $1$ and $U_2$ has been colored with $2$. Let $A_1,A_2$ be the two elements of $\mathcal{A}$ incident to $U$. Choose a shortest path $v_0,\ldots,v_k$ from $A_1$ to $A_2$ in $G[U\cup A_1\cup A_2]$. We can choose distinct vertices $a,d\in U_1$ and $b\in U_2$ among $v_2,\ldots,v_{k-2}$, which are not incident to $W$. Then we adjust their colors, denoting the resulting coloring by $\widetilde c$. Consider the following three possible colorings of $(a,b,d)$:
\begin{enumerate}
    \item $(1,2,1)$, which contributes $1+1-2\equiv0\pmod{3}$ to $\sum_{x\in U_1}\widetilde c(x)-\sum_{x\in U_2}\widetilde c(x)$;
    \item $(3,2,3)$, which contributes $3+3-2\equiv1\pmod{3}$ to $\sum_{x\in U_1}\widetilde c(x)-\sum_{x\in U_2}\widetilde c(x)$;
    \item $(3,2,1)$, which contributes $3+1-2\equiv2\pmod{3}$ to $\sum_{x\in U_1}\widetilde c(x)-\sum_{x\in U_2}\widetilde c(x)$.
\end{enumerate}
Only $a$ and $d$ may change color. Since they belong to the same bipartition class and have no neighbors in $W$, each choice preserves properness, including when these changes are made simultaneously for all $U\in\mathcal{S}_2$.

Note that, for each $U\in\mathcal{S}_2$, all boundary edges have already been assigned weights by $w$, and we keep these weights fixed. For $x\in U$, let $\beta_U(x)=\sum_{\substack{y\notin U\\xy\in E}}w(xy)$. We subtract this boundary contribution from the desired color and choose the unique one of the three colorings for which $\sum_{x\in U_1}(\widetilde c(x)-\beta_U(x))\equiv\sum_{x\in U_2}(\widetilde c(x)-\beta_U(x))\pmod{3}$. By Lemma~\ref{eq}, we can then reassign the internal edges of $U$ so that $s_w(x)\equiv\widetilde c(x)\pmod{3}$ for every $x\in U$. All these finite choices can be made Borelly.

These changes affect no edge incident to $D$ or $W$, so the equations on $D$ are still satisfied. For $x\in W\setminus D$, no incident edge meets $D$, so all its incident edges retain weight $3$, giving $s_w(x)\equiv0\equiv\widetilde c(x)\pmod{3}$. Thus the resulting Borel weighting satisfies $s_w(x)\equiv\widetilde c(x)\pmod{3}$ for every $x\in X$. Since $\widetilde c$ is a proper coloring, this is a Borel vertex-coloring $3$-edge weighting.
\end{proof}

As an example of two-ended bipartite graph, we have the following corollary.
\begin{corollary}\label{line}
Suppose $G=(X, E)$ is a locally finite Borel graph. If each component of $G$ is isomorphic to a line, then $G$ admits a Borel vertex-coloring $3$-edge weighting.
\end{corollary}

This completes the proof of corollary \ref{Zd}.

\subsection{Graphs with Borel selections of ends}
In this subsection, we show more connection between end theory and Borel vertex-coloring edge weighting. Recall that, we say a Borel graph $G=(X,E)$ is \emph{generated by a Borel function $f$} if 
\[
    xEy\Longleftrightarrow x\ne y, (f(x)=y \mbox{ or } f(y)=x).
\]

Then we has the following observation.

\begin{lemma}\label{select}
Let $G=(X,E)$ be a locally finite acyclic Borel graph. Then the following are equivalent:
\begin{enumerate}
    \item $G$ is generated by an aperiodic Borel map $f:X\longrightarrow X$.
    \item There is a Borel selection of one end from each connected component of $G$.
\end{enumerate}
\end{lemma}
\begin{proof}
$(1)\implies(2)$. For each $x\in X$, since $f$ is aperiodic, the forward orbit $\rho_x=(f^n(x))_{n\in\mathbb{N}}$ is a ray. Let $A$ be the set consisting of all rays $\alpha\in[G]^\infty$ such that $f(\alpha(n))=\alpha(n+1)$ for all sufficiently large $n$. Then $A$ is $E_T$-invariant. Note that if $x$ and $y$ are in the same component of $G$, then $\rho_x,\rho_y$ are tail equivalent, so $A\cap[G[C]]^\infty$ is exactly one $E_T$-equivalence class for every connected component $C$.

$(2)\implies(1)$. Let $A\subseteq[G]^\infty$ be a Borel selection of one end from each component. Since each component is a tree, for every $x\in X$ there is a unique ray $\rho_x\in A$ starting at $x$. Note that $x\longmapsto \rho_x$ is Borel by Lusin--Novikov uniformization~\cite[Theorem~18.10]{Kec95}. Define $f(x)=\rho_x(1)$. Then $f$ is Borel and $G_f\subseteq G$. Conversely, let $xy$ be an edge of $G$. If $f(x)\neq y$, then $y$ does not occur in $\rho_x$. Hence $(y,\rho_x(0),\rho_x(1),\ldots)$ is a ray tail equivalent to $\rho_x$ and therefore belongs to $A$. By the uniqueness of the ray in $A$ starting at $y$, this ray equals $\rho_y$. Thus $f(y)=x$ and $G=G_f$. Finally, $f$ is aperiodic because $\rho_x$ is injective.
\end{proof}

Based on the lemma above, we can show the following stronger theorem which implies that function graphs admit a Borel vertex-coloring $3$-edge weighting.

\begin{theorem}\label{selectend}
Suppose $G=(X, E)$ is a locally finite acyclic Borel graph. If there is a Borel way to select a nonempty finite set of ends from each component of $G$, then there is a Borel vertex-coloring $3$-edge weighting.
\end{theorem}

\begin{proof}
By~\cite[Theorem~6.2]{Mil09}, if we can select at least three (finitely many) ends from each component, then the graph is smooth and we can apply Proposition~\ref{smooth}. So it suffices to prove that there is a Borel way to select one or two ends from each component of $G$. We may also assume that no component of $G$ is isomorphic to a line, since those components are handled by Corollary~\ref{line}.

Let $H$ be a Borel spanning subgraph of $G$ defined by:
\begin{center}
$xHy\iff xGy$ and $(\deg(x)\le 3 \deg(y) $ and $\deg(y)\le 3\deg(x))$.
\end{center}
We first show that $H$ is Borel 3-colorable. By Lemma~\ref{select}, we may suppose that $G$ is a Borel graph induced by a finite-to-one aperiodic Borel map $f:X\longrightarrow X$. Note that any subgraph of this graph is still a function graph. For instance, let $f_H:X\longrightarrow X$, where $f_H(x)=x$ if $x$ is not incident to $f(x)$ in $H$; otherwise, let $f_H(x)=f(x)$. Then $H$ is a Borel graph generated by $f_H$. By~\cite[Theorem~5.1]{KST99}, it suffices to show that $H$ admits a finite Borel vertex coloring. For each $i\in{\mathbb{N}}$, let
\begin{center}
$ X_i=\{x: 4^i \le \deg (x) < 4^{i+1}\}$ and $A=\bigcup_{i}X_{2i}$, $B=\bigcup_{i}X_{2i+1}$.
\end{center}
$H[X_i]$ is a Borel graph of bounded degree, so by~\cite[Proposition~4.6]{KST99}, it has a finite Borel coloring. Hence $H[X_i]$ is Borel $3$-colorable. We can color $H[A]$ with $1,2,3$ and $H[B]$ with $4,5,6$, because if $x\in{X_i}$ and $y\in{X_{i+2}}$, then they are not incident in $H$ by definition. This gives a Borel $6$-coloring of $H$. Again by~\cite[Theorem~5.1]{KST99}, there is a Borel coloring $c:X\longrightarrow \{1,2,3\}$ of $H$.

We define a one-ended partial function contained in $G$. The proof is similar to that of~\cite[Theorem~3.7]{CMT16}. Let $A$ be the set of all leaves in $G$ and $S=[A]_{E_G}$ be its saturation. Applying~\cite[Proposition~3.1]{CMT16} to $A$, we obtain a one-ended Borel partial function $q:S\setminus A \longrightarrow S$. Then we consider the set $Y=X\setminus S$. Because $\deg(x)\ge 2$ for each $x\in{Y}$, the preimage $f^{-1}(x)$ is a nonempty finite set. By the Lusin--Novikov theorem~\cite[Theorem~18.10]{Kec95}, we can define a Borel map $g:Y\longrightarrow Y$ such that $f(g(x))=x.$ Note that $g$ is aperiodic and injective; otherwise, $f$ would be periodic. Consider the $G$-invariant Borel set $Z\subseteq Y$ consisting of the vertices in those components of $G[Y]$ on which $g$ is not one-ended. Then $g$ generates a unique bi-infinite line in each component of $G[Z]$. Let $L$ be the Borel set of vertices on these lines. On $Y\setminus Z$, the function $g$ is one-ended. Define $B_n=\{x\in Z:d_G(x,L)\geq n\}$. Note that $(B_n)_n$ decreases to the empty set, and $B_{n+1}$ meets every component of $G[B_n]$: no component of $G[Z]$ is a line, and every vertex of $G[Z]$ has degree at least two. By~\cite[Lemma~3.2]{CMT16}, we obtain a one-ended function $r:Z\longrightarrow Z$.

Finally, we combine these one-ended functions. We define $p:X\setminus A\longrightarrow X$ as follows:
\begin{center}
$ p(x)= \begin{cases}
q(x),&x\in S\setminus A,\\
g(x),&x\in Y\setminus Z, \\
r(x),&x\in Z.
\end{cases} $
\end{center}
$p$ is a one-ended Borel function. Applying Lemma~\ref{partial} to $p$ and $c|_{X\setminus A}$, we obtain a Borel edge weighting $w:E\longrightarrow\{1,2,3\}$ such that $s_w(x)\equiv c(x)\pmod{3}$ for every $x\notin A$. We verify that $w$ is a vertex-coloring $3$-edge weighting. Let $xy\in E$. If $x\in A$, then $y$ has another neighbor because $G$ has no isolated edges. Thus $s_w(x)=w(xy)<s_w(y)$. Now assume that $x,y\notin A$. If $xy$ is an edge of $H$, then $c(x)\neq c(y)$, so $s_w(x)\not\equiv s_w(y)\pmod{3}$. Otherwise, without loss of generality, we have $\deg(x)>3\deg(y)$ by the definition of $H$. Hence $s_w(x)\geq\deg(x)>3\deg(y)\geq s_w(y)$.

We now consider the case where two ends are selected from each component. Let $L$ be the Borel set consisting of the unique lines joining the two selected ends in each component. Choose a Borel proper coloring $b:L\longrightarrow\{1,2,3\}$ of $G[L]$ by~\cite[Proposition~4.6]{KST99}. For each $x\in X$, let $\pi(x)$ be its unique nearest vertex in $L$, and define the Borel map $c:X\longrightarrow\{1,2,3\}$ by $c(x)\equiv b(\pi(x))+d_G(x,L)\pmod{3}$. It is easy to see that $c$ is a proper coloring of $G$. The rest is the same. Again, we let $A$ be the set of leaves of $G$, let $S$ be its saturation in $G$, and put $Y=X\setminus S$. As in the preceding argument involving lines, there is a one-ended Borel function $r:Y\longrightarrow Y$ contained in $G[Y]$. By Lemma~\ref{partial}, there is a Borel weighting $w:E\longrightarrow\{1,2,3\}$ satisfying $s_w(x)\equiv c(x)\pmod{3}$ for every $x\notin A$. Again, by the preceding degree argument, $w$ satisfies our requirements.
\end{proof}

\begin{remark}
    Since every locally finite graph generated by a Borel function has at most one cycle in each connected component, Proposition \ref{smooth} and Theorem \ref{selectend} imply that every  graph generated by a single Borel function admits a Borel vertex-coloring $3$-edge weighting.
\end{remark}

The next corollary uses the topology on the space of ends. For more details, please see~\cite[Sections~1 and~4]{HM09}.

\begin{corollary}\label{infiniteends}
Suppose $G=(X,E)$ is a locally finite acyclic Borel graph. If every connected component of $G$ has at most countably many ends, then $G$ admits a Borel vertex-coloring $3$-edge weighting.
\end{corollary}
\begin{proof}
Note that selecting the whole end space gives a Borel selection of nonempty closed countable sets of ends. By~\cite[Theorem~7.3]{HM09}, there is a Borel selection of one end or a line from each such component. Since $G$ is acyclic, each selected line determines two distinct ends, which gives a Borel selection of two ends. So Theorem~\ref{selectend} applies.
\end{proof}

\subsection{Applications to \texorpdfstring{$\mu$}{mu}-measurable and Baire measurable edge weightings}

Finally, we establish the following two results, which show that the analogue of Theorem \ref{kcolor} holds in the p.m.p. and Baire measurable settings. Let $\mu$ be a Borel probability measure on $X$, and let $\tau$ be a compatible Polish topology. A coloring is \textit{$\mu$-measurable}, respectively \textit{Baire measurable}, if it is Borel outside a Borel $\mu$-null, respectively $\tau$-meager, subset of $X$. Following the edge-coloring convention in~\cite[Section~3]{KM20}, an edge weighting $w:E\longrightarrow\{1,\ldots,k\}$ is \textit{$\mu$-measurable}, respectively \textit{Baire measurable}, if $w\upharpoonright E(G[Y])$ is Borel for some $G$-invariant conull, respectively comeager, Borel set $Y\subseteq X$.

We call $G$ \textit{probability measure-preserving} (\textit{p.m.p.}) if every Borel bijection between Borel subsets of $X$ whose graph is contained in $E_G$ preserves $\mu$. We call $\tau$ \textit{$G$-quasi-invariant} if $[N]_{E_G}$ is meager whenever $N\subseteq X$ is Borel and meager. Before we state and prove the main theorems for this subsection, the following technical lemmas dealing with non-bipartite case are needed.

\begin{lemma}\label{cycle}
Suppose $G=(X,E)$ is a graph. For any given odd cycle $(x_0,\dots,x_{2n},x_0)$ and any map $c: \{x_0,\dots,x_{2n},x_0\}\longrightarrow \{1,2,3\}$, there exists an edge weighting $w:E\longrightarrow \{1,2,3\}$ such that $s_w(x_i)\equiv c(x_i) \pmod{3}$ for every $i$.
\end{lemma}
\begin{proof}
We proceed by induction. First, we deal with the triangle $(x_0,x_1,x_2,x_0)$ and a map $c:\{x_0,x_1,x_2\}\longrightarrow \{1,2,3\}$. Suppose $w:E\longrightarrow \{1,2,3\}$ is a map. Our requirements are:
\begin{enumerate}
    \item $w(x_0x_1)+w(x_2x_0)\equiv c(x_0) \pmod{3}$.
    \item $w(x_0x_1)+w(x_1x_2)\equiv c(x_1) \pmod{3}$.
    \item $w(x_1x_2)+w(x_2x_0)\equiv c(x_2) \pmod{3}$.
\end{enumerate}
Note that the equation $2w(x_0x_1)\equiv c(x_0)+c(x_1)-c(x_2) \pmod{3}$ always has a solution for $w(x_0x_1)$. We choose such a value for $w(x_0x_1)$ and calculate the values of $w(x_1x_2)$ and $w(x_2x_0)$ so that our requirements are satisfied.

Suppose the $2n-1$ case has been proved. We consider the cycle $(x_0,\dots,x_{2n},x_0)$ and $c:\{x_0,\dots,x_{2n}\} \longrightarrow \{1,2,3\}$. We add an edge between $x_0$ and $x_3$, so that $(x_0,x_3,x_4,\dots,x_{2n},x_0)$ is a cycle of length $2n-1$. Consider the map $c':\{x_0,\dots,x_{2n}\}\longrightarrow \{1,2,3\}$, where $c'(x)= c(x)$ if $x\ne x_3$ and $c'(x)= c(x_1)+c(x_3)-c(x_2) \mod 3$ if $x=x_3$. By induction, the values of $w$ on the edges $x_3x_4,\dots,x_{2n}x_0$ have been defined with respect to $c'$. We define the weights of the remaining edges in the following order:
\begin{enumerate}
    \item $w(x_0x_1)\equiv w(x_0x_3) \pmod{3}$.
    \item $w(x_1x_2)\equiv c(x_1)-w(x_0x_1) \pmod{3}$.
    \item $w(x_2x_3)\equiv c(x_2)-w(x_1x_2) \pmod{3}$.
\end{enumerate}
It is easy to check that $w$ satisfies our requirements.
\end{proof}

If the cycle $(x_0,\ldots,x_{2n-1},x_0)$ has even length, then the corresponding equations imply $\sum_{i=0}^{n-1}c(x_{2i})\equiv\sum_{i=0}^{n-1}c(x_{2i+1})\pmod3$, which may false for an arbitrary map $c$.

Recall that a component of a graph is not bipartite if and only if it contains at least one odd cycle. We have the following result.
\begin{lemma}\label{notbipartite}
Suppose $G=(X,E)$ is a locally finite Borel graph. If $G$ is Borel $3$-colorable and each component of $G$ contains at least one odd cycle, then $G$ admits a Borel vertex-coloring $3$-edge weighting.
\end{lemma}
\begin{proof}
Let $[C]^{<\infty}$ be the space of odd cycles of $G$ and consider the intersection graph of it. By the locally finiteness, we take a maximal independent set $M\subseteq{[C]}$ such that any two odd cycles in $M$ are disjoint. Then let
\begin{center}
$C=\bigcup M=\{x:\exists n\ge 1, \exists x_1,\dots,x_{2n}$ distinct vertices of $X$ such that $(x,\dots,x_{2n},x)\in{M}\}$.
\end{center}
It is easy to see that $C$ is a Borel subset of $X$ because $G$ is locally finite. By~\cite[Proposition~3.1]{CMT16}, there is a one-ended Borel partial function $p:X\setminus C\longrightarrow X$. By assumption, we can let $c:X\longrightarrow \{1,2,3\}$ be a Borel proper coloring of $G$. We first apply Lemma~\ref{partial} so that the sum function equals to $c \mod 3$ on $X\setminus C$. For $C$, note that the edges incident to $C$ but not in the cycles have already been assigned weights, so we subtract these contributions from $c$. For instance, if these edges contribute $\beta(x)$ at some $x\in{C}$, then we consider $\hat{c}(x)=c(x)-\beta (x)$ and apply Lemma~\ref{cycle} to $\hat{c}$. Finally, we obtain a Borel vertex-coloring $3$-edge weighting.
\end{proof}

Now we are ready to prove the main results.

\begin{theorem}\label{pmp}
Suppose $G=(X,E)$ is a locally finite Borel graph that is p.m.p. with respect to a Borel probability measure $\mu$. If $G$ has a $\mu$-measurable $3$-coloring, then $G$ admits a $\mu$-measurable vertex-coloring $3$-edge weighting.
\end{theorem}

\begin{proof}
By~\cite[Proposition~2.1(1)]{CMT16}, $G$ has a Borel $3$-coloring on a $\mu$-conull $G$-invariant Borel set $Y$.

We divide $Y$ into two $G$-invariant Borel sets $Y_1,Y_2$, consisting of the vertices in components of $G[Y]$ that are not two-ended and two-ended, respectively. By~\cite[Theorem~2.1]{CGMT21} and $G$ is p.m.p., we may assume that $G[Y_1]$ admits a Borel one-ended spanning subforest. Together with the Borel $3$-coloring, Theorem~\ref{oneend} gives a Borel vertex-coloring $3$-edge weighting on $G[Y_1]$.

For $G[Y_2]$, we apply Proposition~\ref{bipartite} to the bipartite part and Lemma~\ref{notbipartite} to the non-bipartite part. Since these sets are $G$-invariant, we obtain a Borel vertex-coloring $3$-edge weighting on $G[Y_1\cup Y_2]$, where $Y_1\cup Y_2$ is $\mu$-conull in $X$. We extend it to the remaining components by choosing ordinary vertex-coloring $3$-edge weightings, which exist by~\cite[Theorem~2]{Sta23}. The resulting weighting is $\mu$-measurable.
\end{proof}

Next, we turn to the Baire measurable case. To ensure that the topological structures are compatible, we need the following simple but useful result.

\begin{lemma}
Suppose $(X,\tau)$ is a topological space and $E$ is an equivalence relation on $X$. If $\tau$ is $E$-quasi-invariant, then any comeager set $C$ contains an $E$-invariant comeager $G_\delta$ set $C'$.
\end{lemma}
\begin{proof}
It suffices to show that any meager set $M$ is contained in an $E$-invariant meager $F_\sigma$ set $M'$. We inductively define a sequence of set by $E$-quasi-invariance as follow: let $M={M_0}$ and choose $M_{n+1}$ such that $[M_n]_E\subseteq{M_{n+1}}$, where $M_{n+1}$ is meager $F_\sigma$. Now let $M'=\bigcup_n M_n=\bigcup_n [M_{n}]_E$, which is an $E$-invariant meager $F_\sigma$ set.
\end{proof}

We will repeatedly use this lemma in the proof of the following theorem without further mention.

\begin{theorem}\label{baire}
Suppose $G=(X,E)$ is a locally finite Borel graph on a Polish space $(X,\tau)$, and $\tau$ is $G$-quasi-invariant. If $G$ has a Baire measurable $3$-coloring, then $G$ admits a Baire measurable vertex-coloring $3$-edge weighting.
\end{theorem}

\begin{proof}
By~\cite[Proposition~2.1(2)]{CMT16}, together with the fact that $\tau$ is $G$-quasi-invariant, $G$ has a Borel $3$-coloring on a $G$-invariant comeager $G_\delta$ set $Y$.

Similarly, we divide $Y$ into three $G$-invariant Borel sets $Y_1,Y_2,Y_3$, consisting of the vertices in components of $G[Y]$ with one end, two ends, and at least three ends (possibly infinitely many), respectively. By the Baire property, after restricting $Y$ to a $G$-invariant comeager $G_\delta$ subset, we may assume that each $Y_i$ is $G_\delta$ in $Y$, and hence Polish.

On the one-ended part $G[Y_1]$,~\cite[Theorem~1.1(i)]{BPZ24} gives a Borel one-ended spanning tree on a $G$-invariant comeager Borel subset of this part. Together with the Borel $3$-coloring and Theorem~\ref{oneend}, this gives a Borel vertex-coloring $3$-edge weighting on $G[Y_1]$.

On the part $G[Y_3]$ with at least three ends, by~\cite[Theorem~12.1]{KM04}, there is a $G$-invariant comeager $G_\delta$ subset $C$ of this part such that $E_G\upharpoonright C$ is hyperfinite. By~\cite[Lemma~2.4]{Mil09}, there is a Borel spanning subtree $T\subseteq G[C]$ with $E_T=E_G\upharpoonright C$. Every component of $T$ has at least three ends, because a spanning tree will not has less ends than the original graph. By~\cite[Theorem~1.5(2)]{CMT16}, we obtain a Borel one-ended spanning subforest on a $G$-invariant comeager Borel subset of $C$. We then apply Theorem~\ref{oneend} again.

For $G[Y_2]$, we consider two cases: bipartite and non-bipartite. For the bipartite part, we apply Proposition~\ref{bipartite}. For the non-bipartite part we apply Lemma~\ref{notbipartite}.

Since these sets are $G$-invariant, we obtain a Borel vertex-coloring $3$-edge weighting on a comeager set. We can extend it to the whole space as in the preceding proof.
\end{proof}

\section{Further questions}\label{sec:questions}

We conclude this paper with some natural questions concerning the 1--2--3 Conjecture in Borel combinatorics.

\subsection{Borel 3-colorability and 3-edge weightings}

Theorems~\ref{pmp} and~\ref{baire} establish analogues of the classical result in Theorem~\ref{kcolor} for $3$-colorable graphs in the p.m.p. and Baire measurable settings, under their respective hypotheses. It is natural to ask whether Borel $3$-colorability suffices to obtain a Borel vertex-coloring $3$-edge weighting, without the additional structural assumptions used in our Borel results.

\begin{question}
Suppose $G=(X,E)$ is a locally finite Borel graph without isolated edges. If $G$ is Borel $3$-colorable, does $G$ admit a Borel vertex-coloring $3$-edge weighting?
\end{question}

\subsection{Acyclic graphs with continuum many ends}

Let $\mathbb{F}_2=\langle a,b\rangle$, and consider Schreier graphs with respect to the standard generating set $\{a,a^{-1},b,b^{-1}\}$. Every connected component of a Schreier graph arising from a free Borel action of $\mathbb{F}_2$ is a $4$-regular tree and has continuum many ends. Example~\ref{F2} shows that such graphs need not admit a Borel vertex-coloring $2$-edge weighting, whereas Theorem~\ref{regular} guarantees that four weights suffice. This leaves the following question.

\begin{question}
Let $G$ be the Schreier graph induced by a free Borel action of $\mathbb{F}_2$, with respect to the standard generators. Does $G$ admit a Borel vertex-coloring $3$-edge weighting?
\end{question}

In view of Corollary~\ref{infiniteends}, it would be interesting to understand how the structure of the ends affects the existence of Borel vertex-coloring $3$-edge weightings. These $4$-regular trees provide a natural setting in which to investigate this issue.
\subsection{Borel vertex-coloring 2-edge weightings of grid}

The Schreier graphs of free Borel $\mathbb{Z}^d$-actions form an important class of examples in Borel combinatorics. Gao, Jackson, Krohne, and Seward~\cite{GJKS24} developed orthogonal marker decompositions to study Borel colorings, while Gao and Jackson~\cite{GJ26} studied continuous and Borel rectangular partitions. These methods provide a natural starting point for investigating Borel vertex-coloring $2$-edge weightings of grid graphs.

\begin{question}
Let $G$ be the Schreier graph induced by a free Borel action of $\mathbb{Z}^d$ with respect to the standard generators, where $d\geq2$ is an integer. Does $G$ admit a Borel vertex-coloring $2$-edge weighting?
\end{question}

\subsection{The complexity of Borel vertex-coloring 2-edge weightings}

Dudek and Wajc~\cite[Theorem~1.1]{DW11} proved that deciding whether a finite graph admits a vertex-coloring $2$-edge weighting is NP-complete. It is easy to see the definition of Borel vertex-coloring edge weighting is $\boldsymbol{\Sigma}^{1}_{2}$. This suggests asking whether its Borel analogue is $\boldsymbol{\Sigma}^{1}_{2}$-complete, as is the case for certain constraint satisfaction problems; see~\cite[Remark~1.5]{GH26}.

\begin{question}
Is the set of (codes of) locally finite Borel graphs without isolated edges that admit a Borel vertex-coloring $2$-edge weighting $\boldsymbol{\Sigma}^{1}_{2}$-complete?
\end{question}

\section*{Acknowledgements}
The authors would like to thank Longyun Ding and Su Gao for their guidance and support.

\end{document}